\documentclass[french,12pt]{amsart}
\usepackage{babel,graphicx,amssymb,amsmath}
\newtheorem{theoreme}{Th\'eor\`eme}[section]
\newtheorem{corollaire}[theoreme]{Corollaire}
\newtheorem{lemme}[theoreme]{Lemme}
\newtheorem{proposition}[theoreme]{Proposition}
\theoremstyle{definition}
\newtheorem{definition}[theoreme]{D\'efinition}

\theoremstyle{remark}
\newtheorem{remarque}[theoreme]{Remarque}
\newtheorem{exemple}[theoreme]{Exemple}
\numberwithin{equation}{section}

\newcommand{\dbd}{\overline{\partial}\partial}

\begin{document}

\title[Laminations projectives]{Laminations dans les espaces projectifs complexes}
\author{Bertrand Deroin}%
\address{\newline
Max Planck Institut - Leipzig} 
\email{\newline
Bertrand.Deroin@mis.mpg.de}

\thanks{Ce travail a \'et\'e partiellement support\'e par la JSPS et le Fond National Suisse}

\begin{abstract} Dans ce travail, nous \'etendons le Th\'eor\`eme de plongement de 
K.~Kodaira aux vari\'et\'es complexes hermitiennes non compactes 
et aux laminations par vari\'et\'es complexes.  \end{abstract}
\maketitle

\section*{Table des mati\`eres}

1. Introduction.

2. \'Enonc\'e des r\'esultats.

3. Section \`a d\'ecroissance exponentielle. 

4. S\'eries fuchsiennes.

5. Immersion d'une vari\'et\'e complexe non compacte. 

6. Continuit\'e des s\'eries fuchsiennes. 

7. S\'eries fuchsiennes sur une lamination. 

8. Immersion de laminations par vari\'et\'es complexes. 

9. Le cas de la dimension $1$. 

R\'ef\'erences.

\section{Introduction}

Ce travail est une extension aux cas des vari\'et\'es 
hermitiennes non compactes et des laminations par vari\'et\'es complexes 
du Th\'eor\`eme de plongement de K.~Kodaira~\cite{Kodaira1}. 
Ce dernier affirme que parmi les vari\'et\'es complexes \textit{compactes} 
les vari\'et\'es projectives complexes sont celles qui admettent un
fibr\'e en droites holomorphe muni d'une m\'etri\-que hermitienne
dont la courbure est strictement positive. La courbure
d'une m\'etrique hermitienne lisse $|.|$ est d\'efinie par
\[  \Omega = i \dbd \log |s|,\ \ \ \ \ \ i=\sqrt{-1}\]
o\`u $s$ est une section holomorphe locale ne s'annulant pas.

Nos motivations proviennent de l'\'etude qualitative des feuilletages holomorphes singuliers sur les vari\'et\'es 
projectives complexes. Consid\'e\-rons par exemple
une distribution de dimension $d$ sur ${\bf C}^{N+1}$, d\'efinie par des formes holomorphes polynomiales homog\`enes. 
Si elle est int\'egrable, ce qui est toujours le cas si $d=1$, ses vari\'et\'es int\'egrales dessinent un feuilletage
holomorphe singulier invariant par l'action des homoth\'eties et se projette donc en un feuilletage holomorphe singulier 
sur l'espace projectif ${\bf C}P^N$. En fait, tous les feuilletages holomorphes singuliers sur les vari\'et\'es projectives 
sont d\'efinis par des \'equations alg\'ebriques, d'apr\`es le c\'el\`ebre Th\'eor\`eme GAGA de J.~P.~Serre \cite{Serre}. 

La structure g\'eom\'etrique et dynamique de ces feuilletages est tr\`es riche. 
Leurs vari\'et\'es int\'egrales, appel\'ees \textit{feuilles}, 
sont des vari\'e\-t\'es complexes g\'en\'era\-lement non compactes 
immerg\'ees dans la vari\'et\'e ambiante. Leur g\'eom\'e\-trie d'une part, 
et la mani\`ere dont elles se comportent \`a l'infini d'autre part, 
ne sont pas encore bien comprise. En fait, en dehors 
des feuilletages d\'efinis par une \'equation de Riccati, ou de leurs d\'eriv\'es, nous manquons d'exemples 
de feuilletages holomorphes singuliers que nous savons d\'ecrire qualitativement. 

L'\textit{ensemble limite} d'une feuille d'un feuilletage holomorphe est l'ensemble des points d'accumulation des 
lacets tendant \`a l'infini dans la feuille. Il h\'erite d'une structure de 
\textit{lamination par vari\'et\'es complexes} en dehors de l'ensemble singulier. Par exemple, les 
feuilletages de Riccati poss\`edent un unique ensemble limite \textit{exceptionnel}, les autres \'etant des 
feuilles alg\'ebriques. L'ensemble limite est dans ce cas  
le lieu o\`u la dynamique est concentr\'ee. 
  
Il y a de nombreux exemples de laminations par vari\'et\'es 
complexes abstraites, provenant de la th\'eorie de l'it\'eration, des pavages, 
ou de l'arithm\'etique (voir le survol d'\'E.~Ghys~\cite{Ghys1} sur les laminations 
par surfaces de Riemann, et le paragraphe \ref{conjecture}). 
R\'ealiser l'une d'entre elles comme l'ensemble limite d'une feuille d'un
feuilletage holomorphe sur une vari\'et\'e projective complexe donnerait de nouveaux exemples int\'e\-ressants 
de feuilletages holomorphes. 

Nous n'avons malheureusement pas pu r\'ealiser ce ``r\^eve". 
Cependant, l'existence de feuilletages holomorphes sur les vari\'et\'es projectives 
mon\-trent qu'il y a de nombreuses sous-vari\'et\'es complexes et laminations par vari\'et\'es complexes dans les espaces 
projectifs complexes. Ceci est tout \`a fait int\'eressant en soi, et nous nous proposons de les caract\'eriser.  

\section{\'Enonc\'e des r\'esultats}

\subsection{Sous-vari\'et\'es complexes holomorphiquement immerg\'ees dans un espace projectif complexe}
Nous nous int\'eressons d'abord au cas des vari\'et\'es hermitiennes non compactes \`a g\'eom\'etrie born\'ee. 

\begin{definition} \label{localement bilipschitzienne} Soit $(M,g)$ une vari\'et\'e complexe hermitienne. 
Une immersion $\pi :(M,g)\rightarrow
{\bf C}P^N$ est 
dite \textit{localement bilipschitzienne} si il existe un r\'eel $r>0$ et une
constante $K\geq 1$ telle que la restriction de $\pi$ \`a une
boule de rayon $r$ est un plongement $K$-bilipschitzien.\end{definition} 

Cette d\'efinition ne d\'epend pas 
de la m\'etrique utilis\'ee sur ${\bf C}P^N$, puisque ce dernier est compact. 

\begin{theoreme}[Cas d'une vari\'et\'e non compacte] 
Soit $(M,g)$ une vari\'et\'e
hermitienne compl\`ete et $E\rightarrow M$ un fibr\'e en droites
holomorphe admettant une m\'etrique hermitienne \`a g\'eom\'etrie
born\'ee dont la courbure $\omega$ v\'erifie des in\'egalit\'es du
type
\[  \frac{1}{C} g(u) \leq \Omega(u,\sqrt{-1}u)\leq C g(u)\]
pour une constante $C\geq 1$ ne d\'ependant pas du vecteur tangent
$u$. Alors il existe un entier $N$ et une immersion holomorphe localement bilipschitzienne 
$\pi: (M,g)\rightarrow {\bf C}P^N$.\label{immersion d'une variete non compacte}\end{theoreme}

Dans le cas des surfaces riemaniennes orient\'ees, 
qui h\'eritent d'une structure de surface de Riemann hermitienne 
d'apr\`es le Th\'eor\`eme d'Ahl\-fors-Bers \cite{A-B}, nous d\'emon\-trons le r\'esultat suivant.

\begin{theoreme}[Cas des surfaces riemanniennes orient\'ees] 
Une surface riemannienne \`a g\'eom\'etrie born\'ee s'immerge holomorphi\-que\-ment et localement bilipschitziennement 
dans un espace projectif complexe.\label{surfaces orientees} \end{theoreme}

Le Th\'eor\`eme \ref{immersion d'une variete non compacte}, ainsi que la D\'efinition \ref{localement bilipschitzienne}
apparaissent dans le travail \cite{Gromov1} de M. Gromov. La d\'emonstration fait 
usage de la m\'ethode des \textit{s\'eries fuchsiennes},
\'elabor\'ee par H. Poincar\'e \cite{Poincare1}, et qui consiste \`a sommer des sections holomorphes d'un fibr\'e en droites 
holomorphes $E\rightarrow M$ qui ont de bonnes propri\'et\'es de d\'ecroissance \`a l'infini. M. Gromov utilise des sections holomorphes 
de norme $L^1$ finie, mais il nous semble que la convergence des s\'eries fuchsiennes correspondantes n'est \'etablie 
que dans le cas o\`u la vari\'et\'e $M$ est le rev\^etement universel d'une vari\'et\'e projective. Notre contribution 
est l'emploi de sections holomorphes \`a d\'ecroissance exponentielle, qui assure la convergence des s\'eries fuchsiennes associ\'ees
en toute g\'en\'eralit\'e. 
Nous montrons n\'eanmoins la convergence de certaines s\'eries fuchsiennes consid\'er\'ees par M. Gromov, 
mais qui ne suffisent pas pour d\'emontrer \ref{immersion d'une variete non compacte} (voir Remarque \ref{section minimisante L1}).

\subsection{Laminations par vari\'et\'es complexes projectives} 

Une lamination par vari\'et\'es complexes d'un espace topologique $X$ est un atlas complet
$\mathcal L$ d'hom\'eomor\-phismes $\varphi : U\rightarrow B\times T$ d\'efinis sur des ouverts $U$ 
de $X$ et \`a valeurs dans le produit de la boule unit\'e $B$ de ${\bf C}^n$ par un espace topologique $T$,
en sorte que les changements de cartes pr\'eservent la fibration horizontale par boules, et soient biholomorphes 
en restriction aux fibres. 

\begin{definition} Une lamination par
vari\'et\'es complexes ${\mathcal L}$ d'un espace compact $X$ 
est \textit{projective} si
l'ensemble des applications continues $\pi : X\rightarrow {\bf C}P^N$
immergeant holomorphiquement les feuilles s\'eparent
les points de $X$.\end{definition} 

Voici notre r\'esultat principal dans le cadre des laminations par vari\'e\-t\'es complexes. 

\begin{theoreme} [Cas d'une lamination]
Soit $\mathcal L$ une lamination par
vari\'et\'es complexes d'un espace compact n'ayant pas de cycle \'evanouis\-sant.
Si $\mathcal L$ admet un fibr\'e en droites holomorphes de courbure 
strictement positive le long des feuilles, 
alors ${\mathcal L}$ est projective.\label{theoreme de Kodaira lamine}\end{theoreme}

Un \textit{cycle \'evanouissant} est
un lacet contenu dans une feuille et non 
homotope \`a un point dans sa feuille, qui peut \^etre approch\'e dans la topologie uniforme 
par des lacets contenus dans des feuilles proches qui, eux sont homotopes \`a un point
dans leur feuille. Nous ne savons pas d\'emontrer qu'une lamination projective n'a pas de 
cycle \'evanouissant. En revanche les feuilletages holomorphes par courbes d'une vari\'et\'e projective 
n'ont pas de cycle \'evanouissant (voir \cite{Brunella}).

En adaptant la m\'ethode des s\'eries fuchsiennes au cas feuillet\'e, 
\'E. Ghys construit des fonctions
m\'eromorphes sur les laminations hyperboliques \footnote{Une lamination \textit{hyperbolique} est une lamination 
par surfaces de Riemann dont toutes les feuilles sont rev\^etues par le disque unit\'e.} d'un espace compact, 
ayant une transversale totale \cite{Ghys1}. 
Notre m\'ethode permet de contourner l'hypoth\`ese d'existence d'une transversale totale, et de d\'emontrer  
qu'une lamination hyperbolique d'un espace compact a toujours une \textit{multi-transversale totale}.

Nous conjecturons qu'une lamination par vari\'et\'es complexes d'un espace compact de
dimension topologique finie et v\'erifiant les hypoth\`eses du Th\'eor\`eme \ref{theoreme de Kodaira lamine} 
admet un \textit{plongement} holomorphe \`a valeurs dans un espace projectif complexe.
T.~Ohsawa et N.~Sibony \cite{Ohsawa-Sibony1} d\'emontrent ce
fait pour des feuilletages lisses 
par vari\'et\'es complexes de codimension $1$ d'une vari\'et\'e compacte. 
Voici une version symplectique de notre conjecture.

\begin{theoreme} Soit ${\mathcal L}$ une lamination
compacte par vari\'et\'es complexes n'ayant pas de cycle 
\'evanouissant, et $E\rightarrow {\mathcal L}$ un
fibr\'e en droites holomorphe positif. Alors si la dimension topologique de
$X$ est finie, il existe un plongement symplectique $\pi
:X\rightarrow {\bf C}P^N$ pour la structure symplectique
standard de ${\bf C}P^N$.\label{immersion symplectique}\end{theoreme}

R\'ecemment, A. Ibort et D. Mart\^{\i}nez Torres 
ont d\'emontr\'e le r\'esultat analogue pour les feuilletages par 
vari\'et\'es symplectiques de codimension $1$~\cite{I-M}. 
En codimension quelconque et pour des feuilletages par vari\'et\'es symplectiques 
non k\"ahl\'eriennes il reste ouvert. 

Les Th\'eor\`emes \ref{theoreme de Kodaira lamine} et \ref{immersion symplectique} 
sont d\'emontr\'es au paragraphe \ref{immersion projective d'une lamination}. Pour cela nous 
\'etablissons des propri\'et\'es de continuit\'e des s\'eries fuchsiennes aux paragraphes \ref{continuite des series fuchsiennes}
et \ref{series fuchsiennes sur une lamination}.

\subsection{En dimension $1$}

Toute lamination par surfaces riemanniennes orient\'ees est munie d'une structure de lamination par surfaces de Riemann, 
d'apr\`es le Th\'eor\`eme d'Ahlfors-Bers \cite{A-B}. Signalons par ailleurs que dans \cite{Deroin2}, il est d\'emontr\'e 
que l'espace de Teichm\"uller d'une lamination hyperbolique d'un espace compact est de dimension infinie, si elle 
admet une feuille simplement connexe.
Il y a donc de nombreux exemples de laminations par surfaces de Riemann. Le lecteur pourra se r\'ef\'erer au 
survol d'\'E. Ghys \cite{Ghys1}.  

Si toutes les surfaces riemanniennes orient\'ees \`a g\'eom\'etrie born\'ee sont ``projectives", 
ce n'est plus vrai pour des laminations par surfaces riemanniennes orient\'ees. 
Il y a deux obstructions \`a ce fait, essentiellement \'equivalentes. 

La premi\`ere est le fait que l'intersection
d'une lamination projective avec un hyperplan complexe est un
\textit{diviseur}~: c'est 
une fonction non nulle $m: X\rightarrow {\bf N}$ qui est localement donn\'ee par la
multiplicit\'e d'annulation d'une fonction holomorphe non
constante sur chaque feuille. L'existence d'un diviseur 
est une condition topologique non triviale. Par exemple, la composante de Reeb 
du tore plein n'a pas de diviseur passant par sa feuille 
compacte. 
Nous d\'emontrons le r\'esultat suivant. 

\begin{theoreme}\label{obstruction 1} Une lamination par
surfaces de Riemann d'un espace compact qui n'a pas de cycle \'evanouis\-sant 
est projective si et seulement si elle poss\`ede
un diviseur coupant toutes les feuilles.\end{theoreme}

Un \textit{sol\'eno\"{\i}de} est une lamination dont l'espace transverse est totalement discontinu. 

\begin{corollaire} Un sol\'eno\"{\i}de par surfaces de
Riemann d'un espace compact qui n'a pas de cycle \'eva\-nouis\-sant 
est projectif et tendu.\label{solenoide}\end{corollaire} 

La seconde obstruction \`a la projectivit\'e est de nature homologique et a \'et\'e observ\'ee par \'E. Ghys. 
L'ex\-em\-ple le plus simple est une lamination par surfaces de Riemann $\mathcal L$ d'un espace compact $X$ 
ayant une feuille compacte $M$ homologue \`a $0$ dans $X$. Une telle lamination ne peut \^etre projective, 
car une courbe holomorphe compacte d'un espace projectif complexe n'est pas homologue 
\`a $0$~(elle intersecte un hyperplan complexe positivement). 

Dans le cas g\'en\'eral d'une lamination
par surfaces de Riemann d'un espace compact, cette obstruction 
est exprim\'ee en termes de la th\'eorie 
des cycles feuillet\'es
de D.~Sullivan~\cite{Sullivan1}.
Un \textit{cycle feuillet\'e} est un op\'era\-teur
$T:\Omega^{2}({\mathcal L})\rightarrow {\bf R}$ lin\'eaire, ferm\'e
et strictement positif sur les formes strictement positives.
L'int\'egration sur une feuille compacte est un exemple de cycle
feuillet\'e. 
Dans le cas o\`u aucun cycle feuillet\'e de ${\mathcal L}$ n'est
homologue \`a $0$, nous dirons que ${\mathcal L}$ est \textit{tendue}.
Nous d\'emontrons alors le r\'esultat suivant:~\footnote{Notons que souvent une lamination tendue n'a 
pas de cycle \'evanouissant~: par exemple lorsque l'espace total est une vari\'et\'e.}

\begin{theoreme}\label{obstruction 2} Une lamination par surfaces de Riemann d'un espace compact de dimension
topologique finie qui n'a pas de cycle \'evanouissant 
est projective si et seulement si elle est
tendue.\end{theoreme}

Des Th\'eor\`emes \ref{obstruction 1} et \ref{obstruction 2} d\'ecoule un
r\'esultat de nature compl\`etement topologique.

\begin{theoreme} 
Une lamination compacte par surfaces orient\'ees, de dimension
topologique finie et n'ayant pas de cycle \'evanouissant 
est tendue si et seulement si elle admet une
multi-transversale totale.\label{theoreme topologique}\end{theoreme}

S.~Schwartzman a
d\'emontr\'e ce fait pour des feuilletages de dimension $1$
r\'eelle~\cite{Schwartzman1}, et S.~E.~Goodman pour des feuilletages de
codimension $1$~\cite{Goodman1}. Cet \'enonc\'e a un sens pour des
laminations orient\'ees de dimension quelconque, et nous conjecturons qu'il est encore vrai. \\

\textit{Notations.} Dans le texte, $C$ d\'esigne une constante universelle. Nous faisons usage de la m\^eme notation 
pour des constantes a priori diff\'eren\-tes pour simplifier les notations. 

\vspace{0.5cm}

\textit{Remerciements.} 
Ce travail est une partie de ma th\`ese de doctorat~\cite{Deroin1}
effectu\'ee \`a l'\'Ecole Normale Sup\'erieure de Lyon. 
Je remercie chaleureusement \'E.~Ghys pour avoir dirig\'e mes recherches.

\section{Sections \`a d\'ecroissance exponentielle}\label{sections a decroissance exponentielle}

Soit $M$ une vari\'et\'e complexe et $E\rightarrow M$ un fibr\'e 
en droites holomorphe, muni d'une m\'etrique hermitienne $|.|$. 
La \textit{courbure} de $|.|$ est la $(1,1)$-forme d\'efinie par 
\[ \Omega :=i\dbd \log |s|^2\]
o\`u $s$ est une section holomorphe locale de $E$ et $i=\sqrt{-1}$. 
Une $(1,1)$-forme strictement positive en restriction 
\`a toute droite complexe est dite strictement positive.

Il est bien connu (voir par exemple \cite{Demailly}, p. 307 Proposition 12.10) que si $|.|$ est une 
m\'etrique hermitienne de courbure strictement positive, alors 
au voisinage $V_x$ de tout point $x$ de $M$, il existe une section holomorphe 
$s:V_x\rightarrow E$ dont la norme $|s|=e^{\varphi}$
v\'erifie
\[ \varphi = -|z_1|^2-\ldots-|z_n|^2+o(|z|^2),\]
dans un syst\`eme de coordonn\'ees holomorphes $z=(z_1,\ldots,z_n)$ centr\'e en $x$. 
Soit $g$ la m\'etrique \textit{k\"ahl\'erienne} associ\'ee \`a $\Omega$
par la formule
\[ g(u,v)=2\Omega(u,\sqrt{-1}v),\]
ou bien, ce qui revient au m\^eme
\[ g_x= 2(|dz_1|^2+\ldots+|dz_n|^2).\]
\begin{definition}\label{hypothese de geometrie bornee} Une m\'etrique hermitienne 
$|.|$ de $E$ de courbure strictement positive est dite \`a g\'eom\'etrie 
born\'ee si son \textit{rayon} $r(|.|)$ est strictement positif~:
$r(|.|)$ est le supremum des r\'eels $r\geq 0$ tels que 

(i) Le rayon d'injectivit\'e de $(M,g)$ est uniform\'ement 
minor\'e par $r$ et pour tout point $x$ de $M$ il existe un biholomorphisme 
\begin{equation}\label{voisinage bilipschitzien a une boule euclidienne} z:B_g(x,r)\rightarrow U_x\end{equation}
dans un ouvert $U_x$ de ${\bf C}^n$ envoyant $x$ sur $0$ 
et qui est $2$-bilipschitzien, $U_x$ \'etant muni de la m\'etrique euclidienne standard de ${\bf C}^n$.
  
(ii) Pour tout point $x$ de $M$ il existe une section holomorphe 
$s:B_g(x,r)\rightarrow E$ v\'erifiant pour tout $y$ de $B_g(x,r)$:
\begin{equation}\label{sections locales a poids quadratique} e^{-2d_g(x,y)^2} 
\leq |s(y)|\leq e^{-\frac{1}{2}d_g(x,y)^2}.\end{equation}

(iii) La courbure de Ricci de $g$ est born\'ee uniform\'ement sur $M$.
\end{definition}

Le principe de renormalisation suivant est bien connu, et a port\'e 
ses fruits en g\'eom\'etrie symplectique (voir \cite{Tian1,Donaldson}).

\begin{lemme}[Renormalisation] Si $E\rightarrow M$ est un fibr\'e 
en droites holomorphe hermitien de courbure strictement positive et 
\`a g\'eom\'etrie born\'ee, alors les rayons 
des puissances de $|.|$ v\'erifient pour tout entier $k\geq 0$ 
\[  r(|.|^{\otimes k})\geq \sqrt{k} r(|.|).\]
De plus, la courbure de Ricci des m\'etriques $g_k$ induites par $|.|^{\otimes k}$ 
tend uniform\'ement vers $0$ lorsque $k$ tend vers l'infini.
\label{renormalisation}\end{lemme}

\textit{D\'emonstration.} Pour tout entier $k\geq 1$, les m\'etriques 
$|.|^{\otimes k}$ de $E^{\otimes k}$ sont de courbures 
\[  \Omega_k = k\Omega\]
et sont associ\'ees aux m\'etriques k\"ahl\'eriennes $g_k=kg$ sur $M$. En tout point $x$ de $M$, les coordonn\'ees 
\[z_k=\sqrt{k} z\]
exercent un biholomorphisme $2$-bilipschitzien de $B_{g_k}(x,r\sqrt{k})$ dans l'ouvert $\sqrt{k}U_x$
de ${\bf C}^n$, et les sections $s^k:V_x\rightarrow E^k$ v\'erifient 
\[e^{-2d_{g_k}(x,y)^2} \leq |s^k(y)|\leq e^{-\frac{1}{2}d_{g_k}(x,y)^2},\]
pour tout $y$ de $B_{g_k}(x,r\sqrt{k})$. Ainsi le rayon $r(|.|^{\otimes k})$ est au moins sup\'erieur \`a $\sqrt{k} r(|.|)$. 
De plus la courbure de Ricci de la m\'etrique $g_k$ tend uniform\'ement vers $0$ lorsque $k$ tend vers l'infini. 

\vspace{0.2cm}

Ce Lemme motive l'\'etude des m\'etriques hermitiennes de courbure strictement positive et \`a g\'eom\'etrie born\'ee 
ayant un grand rayon et dont la courbure de Ricci est uniform\'ement petite. Dans ces conditions, nous d\'emontrons 
l'existence de sections holomorphes \`a d\'ecroissance exponentielle prolongeant un jet donn\'e en un point.

L'id\'ee est de construire un prolongement holomorphe $h:M\rightarrow E$ dans $L^2_g(|.|_{x,\alpha})$ 
si $\alpha>0$ est assez petit, o\`u la m\'etrique 
$|.|_{x,\alpha}$ est d\'efinie pour tout point $x$ de $M$ et tout r\'eel $\alpha>0$ par
\[ |.|_{x,\alpha}: = e^{\alpha d_g(x,.)} |.|.\]
En vertu de la formule int\'egrale de Cauchy, 
si $|.|$ est une m\'etrique de courbure strictement positive \`a g\'eom\'etrie born\'ee de rayon $r(|.|)=r$, nous avons 
les \textit{in\'egalit\'es de Garding uniformes}~:

\vspace{0.2cm}

\textit{Pour toute section holomorphe} $\tau : B_g(y,r) \rightarrow E$, 
\begin{equation}\label{inegalite de Garding} |\tau(y)| \leq C(r) \sqrt{\int_{B_g(y,r)} |\tau(z)|^2 dv_g(z)},\end{equation}
\textit{o\`u $C(r)$ est une constante ne d\'ependant que de $r$, et que nous choisirons d\'ecroissante de $r$.}

\vspace{0.2cm}

Une section holomorphe $h:M\rightarrow E$ de $L^2_g(|.|_{x,\alpha})$ admet alors la d\'ecroissance
\[ |h(y)| \leq C|h|_{x,\alpha,2} e^{-\alpha d(x,y)},\]
pour tout $y$ de $M$, o\`u $C$ est une constante ne d\'ependant que de $r$. 
L'existence des sections \`a d\'ecroissance exponentielle 
d\'ecoule donc du Lemme suivant. 
\begin{lemme}[Lemme principal]\label{section L2alphax} Il existe des r\'eels $\alpha>0$ et $r_0>0$ tels que si
$E\rightarrow M$ est un fibr\'e en droites holomorphe muni
d'une m\'etrique $|.|$ de courbure strictement positive \`a g\'eom\'etrie born\'ee 
dont le rayon v\'erifie $r(|.|)\geq r_0$ et pour laquelle la courbure de Ricci de $g$ est uniform\'ement 
born\'ee par $1/4$, alors tout $1$-jet $j$ de section holomorphe 
de $M$ dans $E$ en un point $x$ se prolonge en une section holomorphe $h:M\rightarrow E$ de $L^2_g(|.|_{x,\alpha})$ 
de norme inf\'erieure \`a $C|j|$ o\`u $C$ est une constante universelle.\end{lemme}
Sa d\'emonstration est organis\'ee sous forme de paragraphes. 
\subsection{Perturbation}\label{perturbation}
On pourra consulter \cite{Demailly}, p. 422, Theorem 4.5, 
pour une d\'emonstration du r\'esultat suivant, d\^u \`a L.~H\"ormander : 

\vspace{0.2cm}

\textit{Soit $E\rightarrow M$ un fibr\'e en droites au dessus d'une vari\'et\'e
k\"ahl\'erienne \textit{compl\`ete} $(M,g)$, muni d'une m\'etrique $|.|'$ dont la courbure satisfait la positivit\'e
\[ \Omega '-\kappa_{M} \geq \lambda \Omega, \]
$\kappa_{M}$ d\'esignant la courbure de la m\'etrique du fibr\'e canonique $K_M$ de $M$ 
induite par $g$ et $\lambda$ une constante strictement positive.
Si $u$ est une $(0,1)$-forme \`a valeurs dans $E$,
lisse, $L^{2}_g(|.|')$ et $\overline{\partial}$-ferm\'ee
(c'est \`a dire que $\overline{\partial}u=0$), alors il existe une section $v$ de $E$
lisse et $L^{2}_g(|.|')$ qui v\'erifie $\overline{\partial} v=u$ et l'estim\'ee
\begin{equation}\label{estimee de l'inverse du d-barre}  |v|'_{2} \leq \frac{C}{\lambda} |u|'_{2}, \end{equation}
o\`u $C$ est une constante universelle.}

\vspace{0.2cm}

Supposons que l'on ait une section $s$ lisse et $L^{2}_g(|.|')$ qui soit de
surcroit \textit{presque holomorphe}, dans le sens o\`u
$|\overline{\partial} s|'_{2}<<1$. Nous pouvons inverser $u=\overline{\partial} s$
avec les estim\'ees \ref{estimee de l'inverse du d-barre} : nous obtenons une
section $v$ lisse et $L^{2}_g(|.|')$ telle que $\overline{\partial } (s-v)=0$ et
$|v|'_{2}\leq \frac{C}{\lambda} |u|'_{2}$.
La section $h=s-v$ est holomorphe et nous avons
$|s-h|'_{2} <<1 .$
Nous avons donc perturb\'e une section presque-holomorphe en une section holomorphe.

Dans ce paragraphe nous d\'emontrons que l'on peut appliquer les estim\'ees \ref{estimee de l'inverse du 
d-barre} de l'inverse du $\overline{\partial}$
aux m\'etriques $|.|_{x,\alpha}$ avec $\lambda=1/2$ lorsque $\alpha>0$ est 
assez petit (voir Lemme \ref{estimee de l'inverse du d-barre pour les metriques a poids exponentiel}). 
Il nous faut donc convenablement lisser ces m\'etriques 
qui ne sont \`a priori que continues. C'est le but du Lemme suivant.
\begin{lemme}\label{lissage de la metrique} Il existe une constante universelle $C$ telle que les propri\'et\'es 
suivantes soient v\'erifi\'ees. 
Soit $E\rightarrow M$ un fibr\'e en droites holomorphe et 
$|.|$ une m\'etrique de courbure strictement positive de $E$ \`a g\'eom\'etrie
born\'ee telle que $r(|.|)\geq 1$. Pour
toute fonction $1$-lipshitzienne $\psi :(M,g)\rightarrow {\bf R}$, il existe une fonction lisse $\phi :M\rightarrow {\bf R}$
v\'erifiant $|\phi -\psi|_{\infty}\leq 1$ et $|i\overline{\partial}\partial \phi |_{\infty}\leq C$.\end{lemme}
\textit{D\'emonstration.} La fonction $\phi$ est construite
\`a partir d'un noyau $K$
\[  \phi(x)=\int_{M} K(x,y)\psi(y)dv_{g}(y),\]
o\`u $K:M\times M\rightarrow {\bf R}$ est une fonction
v\'erifiant les propri\'et\'es suivantes~: 

\vspace{0.2cm}

(i) $K$ est lisse et positive,

\vspace{0.2cm}

(ii) le support de $K$ est inclu dans $\{ d(x,y)\leq 1\}$

\vspace{0.2cm}

(iii) les diff\'erentielles
des fonctions $K(.,y)$ jusqu'au deuxi\`eme ordre sont born\'ees
par des constantes universelles, 

\vspace{0.2cm}

(iv) en tout point
$x$ de $M$ on a $\int_{M}K(x,y)dv_{g}(y)=1$. 

\vspace{0.2cm}

D\'emontrons d'abord qu'avec ces
propri\'et\'es la fonction $\phi$ v\'erifie le Lemme.
Nous avons d\'ej\`a
\[ \phi(x)-\psi(x)=\int_{M} K(x,y) (\psi(y)-\psi(x))dv_{g}(y),\]
ce qui donne $|\phi(x)-\psi(x)|\leq 1$
d'apr\`es les propri\'et\'es (ii) et (iv). Pour contr\^oler les d\'eriv\'ees
de $\phi$ en un point $x$, observons que pour tout point $x_0$
\[
\phi(x)-\psi(x_{0})= \int_{M} K(x,y)(\psi(y)-\psi(x_{0}))dv_{g}(y),\]
si bien que, diff\'erentiant et prenant $x_0=x$ nous obtenons
\[ i\overline{\partial}\partial \phi _{|x} = \int_{M}
i\overline{\partial}\partial K(.,y)_{|x} (\psi(y)-\psi(x))dv_{g}(y).\]
Nous avons donc la majoration
\[ |\overline{\partial}\partial \phi _{|x}|\leq C
\mathrm{vol}(B_g(x,1)),\]
la constante $C$ \'etant donn\'ee par (iii). Le r\'esultat
r\'esulte du fait que le volume de la boule $B_g(x,1)$ est major\'e par le volume de la boule euclidienne de rayon $2$.

Pour construire le noyau $K$ nous utilisons des \textit{noyaux centr\'es}
en un point $t$ de $M$. \`A partir des coordonn\'ees $z_t$ 
de~\ref{voisinage bilipschitzien a une boule euclidienne} centr\'ees en $t$,
consid\'erons les fonctions $H_{t}: M\times M\rightarrow {\bf R}$ d\'efinies par 
\[ H_t(x,y) = H(z_t(x),z_t(y))\]
o\`u $H$ est une fonction lisse positive des deux variables $z_1\in {\bf C}^n$ et $z_2\in {\bf C}^n$ v\'erifiant
$H(z_1,z_2)=1$ si $z_1$ et $z_2$ sont de norme inf\'erieure \`a $1/8$ et $H(z_1,z_2)=0$ si l'une des deux normes de 
$z_1$ ou $z_2$ est sup\'erieure \`a $1/4$. Comme l'image $U_t$ de $B_g(t,1)$ par $z_t$ contient $B_{eucl}(0,1/4)$, la fonction 
$H_t:M\times M\rightarrow {\bf R}$ prolong\'ee en dehors de $B_g(t,1)\times B_g(t,1)$ par $0$ est lisse et positive
et v\'erifie les propri\'et\'es suivantes. 
Si $d(x,t)\leq 1/16$ et $d(y,t)\leq 1/16$,
on a $H_{t}(x,y)=1$. Si $d(x,t)\leq 1/4$ ou $d(y,t)\leq 1/4$,
on a $H_{t}(x,y)=1$. A $y$ fix\'e, $i\overline{\partial} \partial H_{t}(.,y)$
est born\'e uniform\'ement par une constante universelle.

Une partie $1/32$-s\'epar\'ee maximale pour l'inclusion
et aussi $1/32$-dense. Prenons en une $T\subset M$. La somme
\[  H=\sum_{t\in T} H_{t}\]
est localement finie car la courbure de Ricci est born\'ee. 
Plus pr\'ecis\'ement, pour tout couple $(x,y)$ de points
de $M$, il n'y a au plus que $|T\cap B_g(x,1/2)|$ points 
$t$ de $T$ pour lesquels $H_t(x,y)\neq 0$. Mais, puisque $T$ est $1/32$-s\'epar\'ee nous avons
\[ |T\cap B_g(x,1/2)| \leq \mathrm{vol}(B_g(x,1/2 +1/32)) / \mathrm{vol}(B_g(x,1/32))\leq 68^{2n}\]
en comparant le volume des boules $B_g$ avec celles des boules euclidiennes de ${\bf C}^n$. 
Ainsi $H$ est une fonction lisse positive sur $M\times M$, dont le support est contenu
dans $\{ d(x,y)\leq 1/2\}$ et dont les d\'eriv\'ees jusqu'au 
deuxi\`eme ordre (y compris $0$) sont born\'ees par des constantes universelles.
D'autre part, si $x$ et $y$ sont deux points de $M$ 
s\'epar\'es d'une distance inf\'erieure \`a $1/32$, il existe un point $t$ de $T$ 
tel que $d(x,t)\leq 1/16$ et $d(y,t)\leq 1/16$. C'est donc que $H(x,y)\geq H_t(x,y)=1$. 
Nous pouvons donc poser
\[ K(x,y)= \frac{H(x,y)}{\int_{M}H(x,y)dv_{g}(y)}.\]
C'est le noyau que nous voulons. Les propri\'et\'es (i), (ii), (iv) sont claires.
Pour d\'emontrer la propri\'et\'e (iii), il suffit de remarquer que
l'int\'egrale $$\int_{M}H(x,y)dv_{g}(y)$$ est minor\'ee uniform\'ement
par une constante strictement positive, et que son
$i\overline{\partial} \partial$ est born\'e par une constante universelle. Le Lemme \ref{lissage de la metrique} est d\'emontr\'e. 

\begin{lemme}\label{estimee de l'inverse du d-barre pour les metriques a poids exponentiel}
Il existe $\alpha>0$ tel que si $|.|$ est une m\'etrique de $E$ de courbure strictement positive \`a g\'eom\'etrie born\'ee 
de rayon $r(|.|)\geq 1$ et pour laquelle la courbure de
Ricci de $g$ est major\'ee par $1/4$, alors pour tout $x$ de $M$ l'estim\'ee suivante est v\'erifi\'ee. 
Si $u$ est une $(0,1)$-forme \`a valeurs dans $E$,
lisse, $L^{2}_{g}(|.|_{x,\alpha})$ et $\overline{\partial}$-ferm\'ee, alors il existe une section $v$ de $E$
lisse et $L^{2}_{g}(|.|_{x,\alpha})$ qui v\'erifie $\overline{\partial} v=u$ et l'estim\'ee
\[  |v|_{x,\alpha,2} \leq C |u|_{x,\alpha,2}, \]
o\`u $C$ est une constante universelle.\end{lemme}
\textit{D\'emonstration.} Consid\'erons la fonction $\psi=d(x,.)$. C'est une fonction 
$1$-lip\-shitzien\-ne. Le Lemme \ref{lissage de la metrique} nous donne une 
fonction $\phi$ v\'erifiant 
$|\phi -d(x,.)|_{\infty}\leq 1$ et 
$|i\overline{\partial}\partial \phi|_{\infty}\leq C$. 
Regardons alors la m\'etrique de $E$ d\'efinie par
\[|.|' = e^{\alpha \phi }|.|.\]
D'une part elle est uniform\'ement \'equivalente \`a $|.|_{x,\alpha}$, avec l'estim\'ee 
$$e^{-\alpha}|.|_{x,\alpha}\leq |.|'\leq e^{\alpha}|.|_{x,\alpha}.$$
D'autre part elle est lisse de courbure
\[ \Omega'=\Omega +\alpha i\overline{\partial}\partial \phi.\]
Si $\alpha >0$ est choisi assez petit en sorte que $\alpha C\leq 1/4$, 
o\`u $C$ est la constante universelle donn\'ee par le Lemme 
\ref{lissage de la metrique}, et si la courbure de Ricci de $g$ est major\'ee par $1/4$,
alors nous avons la positivit\'e de $|.|'$:
\[ \Omega ' - \kappa_{M} \geq  \frac{1}{2} \Omega. \]
Nous pouvons donc utiliser les estim\'ees de H\"ormander \ref{estimee de l'inverse du d-barre} sur $E$
muni des normes $|.|_{x,\alpha}$, avec la constante $\lambda =1/2$ et une constante $C$ universelle.

\subsection{Sections presque-holomorphes} Soit $E\rightarrow M$ un 
fibr\'e en droi\-tes holomorphe et $|.|$ une m\'etrique de courbure strictement
positive \`a g\'eom\'etrie born\'ee. Le Lemme suivant 
est d\^u \`a G.~Tian~\cite{Tian1}.
\begin{lemme}[Tian] \label{lemme de Tian} Soit $x$ un point de $M$ et 
$j$ un $1$-jet de section holomorphe de $M$ dans $E$ en $x$. 
Il existe une section lisse $\overline{s} :M\rightarrow E$ \`a 
support compact, holomorphe sur $B_g(x,r/3)$, passant par $j$, 
telle que 
\[ |\overline s|_{x,\alpha,2} \leq C|j|\ \ \ \ \mathrm{et}
\ \ \ \ |\overline{\partial} \overline{s} |_{x,\alpha,2}\leq \frac{C|j|}{r},\]
o\`u $C$ est une constante universelle.\end{lemme} 
 \textit{D\'emonstration.} Soit $s:B_g(x,r)\rightarrow E$ 
 la section d\'efinie en \ref{sections locales a poids quadratique}, et soit $P$ un
polyn\^ome lin\'eaire des coordonn\'ees $z$ de \ref{voisinage 
bilipschitzien a une boule euclidienne} centr\'ees en $x$ tel que $j=J_1(Ps)(x)$. 
La section $\overline{s}$ est 
d\'efinie par $\overline{s} =\varphi s$, o\`u $\varphi:B_g(x,r)\rightarrow {\bf R}$ 
est une fonction lisse, identiquement \'egale \`a $1$ sur 
$B_g(x,r/3)$, nulle \`a l'ext\'erieur de $B_g(x,2r/3)$. 
De plus nous demandons que $\varphi$ prenne ses valeurs entre $0$ et $1$ et que 
\[ |d\varphi|_g \leq C/r,\]
o\`u $C$ est une constante universelle. Nous pouvons par exemple construire 
$\varphi$ en prenant une fonction $\psi$ ayant les m\^emes
propri\'et\'es sur la boule euclidienne $B_{\mathrm eucl}(0,1/2)$ de ${\bf C}^n$ et en consid\'erant 
$\varphi (y) =\psi (z(y)/2r)$. Nous avons alors 
\[ |\overline{\partial} \overline s|_{x,\alpha,2}=
|\overline{\partial}(\varphi) Ps|_{x,\alpha,2}\leq \frac{C}{r} |Ps|_{x,\alpha,2}.\]
Reste \`a majorer la norme de $Ps$ sur la boule $B_g(x,r)$. Il existe 
une constante universelle $C$ telle que pour tout $y$ de $B_g(x,r)$  
\[ |P(y)| \leq  C|j|(1+d(x,y)).\]
Nous en d\'eduisons donc 
\[  |Ps|_{x,\alpha,2} \leq C|j|\sqrt{\int_{B_g(x,r)}(1+d(x,y))^2e^{\alpha d(x,y)-d(x,y)^2}dv_g(y).}\]
Cette derni\`ere int\'egrale est major\'ee par l'int\'egrale convergente sur ${\bf  C}^n$
\[ 2^{2n} \int_{{\bf C}^n} (1+2|z|)^2 e^{2\alpha |z|- \frac{|z|^2}{2}}dv_{eucl}(z).\] 
Le Lemme est d\'emontr\'e. 

\subsection{D\'emonstration du Lemme principal~\ref{section L2alphax}.}\label{demonstration de la proposition}
Soit $|.|$ une m\'etri\-que de $E$ de courbure strictement positive 
et \`a g\'eom\'etrie born\'ee pour laquelle la courbure de Ricci de $g$ 
est major\'ee par $1/4$, et dont le rayon $r(|.|)$ est sup\'erieur \`a~$1$. 
Soit $x$ un point de $M$ et $j$ un $1$-jet en $x$ de 
section holomorphe de $M$ dans $E$. Appliquons le Lemme de perturbation 
\ref{estimee de l'inverse du d-barre pour les metriques a poids exponentiel}
\`a la section $\overline s$ construite au Lemme \ref{lemme de Tian}. 
Nous obtenons une section holomorphe $H:M\rightarrow E$ de $L^2_g(|.|_{x,\alpha})$ 
v\'erifiant 
\[ |H-\overline s|_{x,\alpha,2}\leq \frac{C|j|}{r}.\]
La section $H-\overline s$ est holomorphe sur la boule $B_g(x,r/3)$, 
donc les in\'egalit\'es de Garding \ref{inegalite de Garding} donnent
\[ |J_1 H -j|\leq \frac{C|j|}{r}\]
pour une constante universelle $C$. Nous choisissons $r$ en sorte 
que $C/r=1/2$, o\`u $C$ est la constante apparaissant dans l'in\'egalit\'e 
pr\'ec\'edente. Puisque d'apr\`es le Lemme de Tian nous avons 
$|\overline s|_{x,\alpha,2} \leq C|j|$, nous en d\'eduisons l'estim\'ee~:
\[ |H|_{x,\alpha,2}\leq C|j|,\]
pour une constante universelle $C$. Ces in\'egalit\'es sont le premier pas d'un 
processus d'approximations successives. D\'efinissons $h_1=H(j)$ et par 
r\'ecurrence $h_{q+1}= H(j-J_1(h_1+\ldots+h_q))$ pour $q\geq 2$. 
Nous avons alors pour tout $q\geq 1$, 
$$|J_1 (h_1+\ldots +h_q) -j|\leq (1/2)^q |j|,\ \ \ |h_{q+1}|_{x,\alpha,2}\leq C(1/2)^q |j|.$$ 
La s\'erie 
$\sum_q h_q$ converge donc vers une section holomorphe $h:M\rightarrow E$ de 
$L^2_g(|.|_{x,\alpha,2})$ passant par $j$ et telle que 
\[ |h|_{x,\alpha,2}\leq C|j|,\]
o\`u $C$ est une constante universelle. Le Lemme~\ref{section L2alphax} est d\'emontr\'e.

\section{S\'eries fuchsiennes}\label{series fuchsiennes}
Soit $E\rightarrow M$ un fibr\'e en droites holomorphe muni d'une m\'etri\-que $|.|$ de courbure strictement positive et 
\`a g\'eom\'etrie born\'ee.  
Soit $g$ la m\'etrique k\"ahl\'erienne associ\'ee \`a la forme de courbure de $|.|$. 
Nous supposons que~:

\vspace{0.2cm}

\ (i) la courbure de Ricci de $g$ est uniform\'ement major\'ee par $1/4$. 

\vspace{0.2cm}

(ii) le rayon $r(|.|)$ de $|.|$ est 
sup\'erieur au r\'eel $r_0$. 

\vspace{0.2cm}

Dans ces conditions, nous avons montr\'e au Lemme~\ref{section L2alphax} que tout $1$-jet $j_x$ 
de section holomorphe de $M$ dans $E$ en un point $x$ se prolonge en une section holomorphe de 
$L^2(|.|_{x,\alpha})$ de norme inf\'erieure \`a $C|j_x|$, o\`u $C>0$ est une constante universelle. 
Soit $m(j_x):M\rightarrow E$ celle de norme minimale.

\begin{definition}
Soit $\delta $ un r\'eel strictement positif. Une partie $T\subset M$ est \textit{$\delta$-s\'epar\'ee} si 
deux points distincts de $T$ sont s\'epar\'es d'une distance sup\'erieure 
\`a $\delta$. Soit $T\subset M$ une partie $\delta$-s\'epar\'ee et $j=\{j_t\}_{t\in T}$ une famille 
de $1$-jets de section holomorphe de $M$ dans $E$ d\'efinis aux points de $T$. La s\'erie 
\[\sigma(j):= \sum_{t\in T} m(j_t) \]
est la \textit{s\'erie fuchsienne} associ\'ee \`a $j$. \end{definition}

\begin{lemme}[Convergence des s\'eries fuchsiennes] 
Il existe une constan\-te $c_0>0$ telle que si la courbure de Ricci de $g$ 
est born\'ee uniform\'ement par $c_0$ alors les s\'eries fuchsiennes
$\sigma(j)$ convergent uniform\'ement sur tout compact de $M$ vers une section holomorphe de $E$ v\'erifiant 
\[ |\sigma(j)|_{\infty,M}\leq C(\delta) |j|_{\infty,T},\]
pour toute partie $\delta$-s\'epar\'ee $T$ et toute famille born\'ee $j$ de $1$-jets d\'efinis sur $T$. De plus, nous avons 
\[ |J_1\sigma(j) -j|_{\infty,T} \leq D(\delta) |j|_{\infty,T},\]
o\`u $D$ tend vers $0$ lorsque $\delta$ tend vers 
l'infini.~\footnote{$|.|_{\infty,A}$ est la norme uniforme en restriction \`a la partie $A$.}
\label{convergence des series fuchsiennes}\end{lemme}
\textit{D\'emonstration.} Remarquons que, pour tout point $y$ de $M$ nous avons les in\'e\-galit\'es 
\[ \sum_{t\in T} e^{-\alpha d(t,y)} \leq C(\delta) \sum_{t\in T} \int_{B_g(t,\delta/2)} e^{-\alpha d(z,y)}dv_g(z)\]
\[\leq C(\delta) \int_M e^{-\alpha d(z,y)}dv_g(z),\]
o\`u $C(\delta) = e^{\alpha \delta}/\nu(\delta/2)$ et 
$\nu(\delta)= \inf _{x\in M} \mathrm{vol} (B_g(x,\delta))$. Il existe une constante $c_0>0$
telle que si la courbure de Ricci de $g$ est uniform\'ement 
major\'ee par $c_0$, alors la croissance du volume des boules 
est exponentielle d'exposant inf\'erieur \`a
$\alpha/2$. Plus p\'ecis\'ement il existe une constante universelle $C$ telle que pour tout $r\geq 1$
\[ \mathrm{vol}(B_g(x,r)) \leq C e^{\alpha r/2}.\]
Alors pour tout $y$ de $M$ l'int\'egrale 
\[ \int_M e^{-\alpha d(z,y)}dv_g(z)\]
est born\'ee par une constante universelle. Nous obtenons donc 
\[ \sum_{t\in T} |h(j_t)(y)| \leq C(\delta ) |j|_{\infty,T},\]
ce qui montre d\'ej\`a que $\sigma$ converge uniform\'ement 
sur tout compact de $M$ vers une section holomorphe de $E$ born\'ee par $C(\delta )|j|_{\infty,T}$.

Supposons maintenant que $\delta \geq 2$. 
En vertu de l'in\'egalit\'e de Garding \ref{inegalite de Garding} nous avons pour tout $y$ de $T$, 
\[ |J_1\sigma (y) -j_y|\leq \sum_{t\in T, t\neq y} 
|J_1 h(j_t)(y)|\leq C|j|_{\infty} \sum_{t\in T,t\neq y} e^{-\alpha d(t,y)}.\]
Nous en d\'eduisons donc 
\[ |J_1\sigma (y) -j_y|\leq C|j|_{\infty}\frac{e^{\alpha}}{\nu(1)}\int_{M-B_g(y,\delta -1)} e^{-\alpha d(z,y)}dv_g(z).\]
Mais les fonctions \[f_{\delta}: y\in M\mapsto \int_{M-B_g(y,\delta -1)} e^{-\alpha d(z,y)}dv_g(z)\]
convergent uniform\'ement vers $0$ lorsque $\delta$ tend vers l'infini. 
Le Lem\-me~\ref{convergence des series fuchsiennes} est d\'emontr\'e.

\begin{proposition}\label{prolongement d'une famille de jets} 
Soit $E\rightarrow M$ un fibr\'e en droites holomorphe muni d'une m\'etrique de courbure strictement 
positive \`a g\'eom\'etrie born\'ee dont le rayon $r(|.|)$ est 
sup\'erieur \`a $r_0$ et pour laquelle la courbure de Ricci de $g$ est uniform\'ement 
major\'ee par $c_0$. Alors il existe $\delta_0>0$ tel que toute famille 
born\'ee de $1$-jets d\'efinie sur une partie $\delta_0$-s\'epar\'ee $T$ de $M$
se prolonge en une section holomorphe $\sigma : M\rightarrow E$ v\'erifiant 
\[ |\sigma|_{\infty,M} \leq C|j|_{\infty,T},\]
o\`u $C$ est une constante universelle. \end{proposition}
\textit{D\'emonstration.} Nous choisissons $\delta _0>0$ tel que $D(\delta_0)\leq 1/2$ et nous 
reproduisons l'argument d'approximation de \ref{demonstration de la proposition}. 
D\'efinissons $\sigma _1= \sigma (j)$, et par r\'ecurrence 
$\sigma_{q+1}=\sigma(j-J_1(\sigma_1+\ldots+\sigma_q))$ pour tout entier $q$ sup\'erieur \`a $1$. 
Nous avons, pour tout $q$ les in\'egalit\'es 
\[  |j-J_1(\sigma_1+\ldots +\sigma_q)|_{\infty,T} \leq |j|_{\infty} (1/2)^q,
\ \ \ |\sigma_q|_{\infty,M} \leq C(\delta_0) |j|_{\infty} (1/2)^{q-1}.\]
Ainsi la s\'erie $\sum_q \sigma_q$ converge uniform\'ement sur $M$ 
vers une section holomorphe $\sigma :M\rightarrow E$ born\'ee par $C|j|_{\infty}$ et 
prolongeant la famille de jets $j$. La constante $C$ est 
universelle car nous pouvons choisir $C$ d\'ecroissante de $\delta$. 

\begin{remarque}\label{section minimisante L1} M.~Gromov propose de fabriquer 
des s\'eries fuchsiennes avec des sections
qui sont seulement dans $L^1_g(|.|)$ (\cite{Gromov1}, p. 387, Corollary 3.3.5). 
Soit $E\rightarrow M$ un fibr\'e en 
droites holomorphe muni d'une m\'etrique de courbure strictement positive 
et \`a g\'eom\'etrie born\'ee.
Si la courbure de Ricci de $g$ est suffisament petite et si le rayon $r(|.|)$ est suffisament grand alors 
par tout point $e$ de $E$ passe une section holomorphe $h$ telle que $|h|_2\leq C|e|$, 
o\`u $C$ est une constante universelle. 
Le produit de deux sections de $L^{2}_{g}(|.|)$ est une section de $E^{\otimes 2}$
de norme $L^{1}$ finie. M.~Gromov consid\'ere
la section holomorphe de norme $L^1$ minimale prolongeant un \'el\'ement de $E^{\otimes 2}$. 
La convergence des s\'eries fuchsiennes
obtenues \`a partir de ces sections est annonc\'ee sans d\'emonstration 
(\cite{Gromov1}, p. 389, Interpolation Theorem 3.3.10).
Nous en donnons une d\'emonstration ici,
avec une variante de la d\'efinition de la section minimisante~: si $e=e_1\otimes e_2$ 
est un \'el\'ement de $E^{\otimes 2}$ nous 
d\'efinissons $m(e) =m(e_1)\otimes m(e_2)$ o\`u $m(e_i)$ est la section holomorphe de $E$ passant par $e_i$ de norme 
$L^2_g(|.|)$ minimisante. 
L'avantage de ces s\'eries fuchsiennes est que le rayon $r(|.|)$ n\'ecessaire pour construire des sections holomorphes $L^2$ de $E$  
est inf\'erieur au $r_0$ de~\ref{section L2alphax} que l'on doit prendre pour construire des 
sections holomorphes \`a d\'ecroissance exponentielle. Par contre nous ne savons 
pas d\'emontrer la convergence des s\'eries form\'ees \`a partir 
des sections minimisantes $L^1$ passant par une famille de $1$-jets.
\begin{lemme} Soit $T\subset M$ une partie s\'epar\'ee, et $\{e_{t}\}_{t\in T}$ une famille born\'ee d'\'el\'ements
$e_{t}$ de la fibre de $E^{\otimes 2}$ au dessus du point $t$. Alors la s\'erie
\[ \sum_{t\in T} m(e_{t}) \]
converge uniform\'ement vers une section holomorphe de $E^{\otimes 2}$. \label{serie fuchsienne3}\end{lemme}
\textit{D\'emonstration.} Ce Lemme vient en fait d'une propri\'et\'e de sym\'etrie des sections minimisantes
$m(e)$ : 

\vspace{0.2cm}

\textit{Propri\'et\'e de sym\'etrie~: soient 
$x,x'$ deux points de $M$ et $e,e'$ deux \'el\'ements de $E$ au dessus de $x$ et $x'$ respectivement.
Alors on a des in\'egalit\'es }
\[\frac{1}{C}|e||m(e')(x)| \leq |e'| |m(e)(x')|\leq C |e||m(e')(x)|,\]
\textit{pour une constante $C$ ne d\'ependant que de $g$.}

\vspace{0.2cm}

D\'emontrons d'abord que cette propri\'et\'e entraine le Lemme~\ref{serie fuchsienne3}.
Remarquons qu'elle est aussi vraie pour les sections $L^{1}$ de
$E^{\otimes 2}$ qui sont les carr\'es de deux sections minimisantes
$L^{2}$ de $E$. Ensuite, prenons un point $x$ dans $M$ et un
\'el\'ement $e$ de $E^{\otimes 2}$ au dessus de $x$ qui est de norme $1$.
Nous avons alors les majorations
\[  
\sum_{t\in T} |m(e_{t})(x)| \leq C (\sup _{t\in T}|e_{t}|)\sum_{t\in T} |m(e)(t)| .\]
Mais d'apr\`es l'in\'egalit\'e de Garding nous avons
\[ \sum_{t\in T} |m(e)(t)| \leq C(\delta) |m(e)|_{1}  .\]
Ainsi $\sum_{t\in T} |m(e_{t})(x)|$ est finie pour tout $x$ et le Lemme~\ref{serie fuchsienne3} en d\'ecoule.

D\'emontrons maintenant la propri\'et\'e de sym\'etrie.
Si $\{s_{\nu}\}_{\nu}$ est une base orthonorm\'ee topologique de
l'espace des sections holomorphes et $L^{2}$ de $E$, les
sections $m(e)$ ont une expression tr\`es simple donn\'ee par
(voir \cite{Bergmann1}) :
\[m(e)= \frac{1}{\sum_{\nu}|s_{\nu}(x)|^{2}}\sum_{\nu} <e,s_{\nu}(x)>s_{\nu},\]
o\`u $x$ est le point sur lequel se projette $e$. Ceci montre en particulier que
\[  |m(e)|_{2}^{2}= \frac{|e|^{2}}{\sum_{\nu}|s_{\nu}(x)|^{2}} ,\]
et que la fonction $x\in M\mapsto \sum_{\nu} |s_{\nu}(x)|^{2}\in {\bf R} _{+}$ 
prend des valeurs comprises entre deux constantes strictement
positives. Fixons deux points $x$ et $x'$ de $M$ et deux \'el\'ements $e$ et $e'$ au dessus de $x$ et $x'$ respectivement.
Nous avons alors
\[ \begin{array} {ccc}
|e'| |m(e)(x')|& = & \frac{|e'|}{\sum_{\nu}|s_{\nu}(x)|^{2}}|\sum_{\nu}<e,s_{\nu}(x)>s_{\nu}(x')|\\
                       & = & \frac{|e|}{\sum_{\nu}|s_{\nu}(x)|^{2}} |\sum_{\nu}<e',s_{\nu}(x')> s_{\nu}(x)|\\
                       & \leq & C|e| |m(e')(x)|, \\
  \end{array} \]
et le Lemme \ref{serie fuchsienne3} est d\'emontr\'e.\end{remarque}

\section{Immersion d'une vari\'et\'e complexe non compacte}
\label{theoreme d'immersion de varietes non compactes}

Dans cette partie nous d\'emontrons les Th\'eor\`emes \ref{immersion d'une variete non compacte} et 
\ref{surfaces orientees} d'immersion de vari\'et\'es
hermitiennes non compactes. 

\subsection{D\'emonstration du Th\'eor\`eme \ref{immersion d'une variete non compacte}}
Soit $E\rightarrow M$ un fibr\'e en droites muni d'une m\'etrique hermitienne $|.|$
de courbure strictement positive et \`a g\'eom\'etrie born\'ee. 
Soit $g$ la m\'etrique k\"ahl\'erienne sur $TM$ associ\'ee \`a la courbure de
$|.|$. Il s'agit de montrer que $(M,g)$ est projective.

 Nous construisons une immersion 
de la forme $ \pi =[\sigma_{0}:\ldots:\sigma_{N}]: M\rightarrow {\bf C}P^N$, 
o\`u les sections $\sigma _i$ sont des sections holomorphes d'une puissance $E^{\otimes k}$ de $E$
pour $k$ assez grand.  Nous avons vu que $r(|.|^{\otimes k})$ est au moins sup\'erieur \`a $\sqrt{k}r(|.|)$, et que 
la m\'etrique $g_k=kg$ est celle associ\'ee \`a la courbure 
$\Omega_k=k\Omega$ de $|.|^{\otimes k}$, si bien que sa courbure de Ricci 
tend uniform\'ement vers $0$ quand $k$ tend vers l'infini. Nous pouvons donc appliquer 
le r\'esultat de la Proposition \ref{prolongement d'une famille de jets} aux fibr\'es $E^{\otimes k}$ et \`a 
leur m\'etrique $|.|^{\otimes k}$~(voir Lemme~\ref{renormalisation}).
Pour simplifier les notations, nous supposons que c'est $E$ lui-m\^eme 
qui v\'erifie les hypoth\`eses de la Proposition~\ref{prolongement d'une famille de jets}.

\begin{lemme}\label{sections holomorphes pour l'immersion} Soient $\delta >0$ le r\'eel donn\'e 
par la Proposition \ref{prolongement d'une famille de jets}.
Il existe un r\'eel $\varepsilon >0$ tel que pour toute partie $\delta$-s\'epar\'ee $T$ de $M$, il existe 
une application m\'eromorphe $\pi_{T}: M\rightarrow {\bf C}P^n$
de la forme $\pi_{T}=[\sigma_{0}:\ldots:\sigma_{n}]$ o\`u les $\sigma_{l}$ 
sont des sections holomorphes de $E$ born\'ees par une constante universelle, telle que
\begin{itemize}
\item $|\sigma_0|\geq 1/2$ sur la r\'eunion $T_{\varepsilon}$ des boules $B_g(t,\varepsilon)$ 
autour des points de $T$, si bien que $\pi_T$ 
est bien d\'efinie sur $T_{\varepsilon}$.
\item sur chaque boule $B_g(t,\varepsilon)$, nous avons $\pi_{T}^{*}\Omega_{FS}\geq C\Omega$, o\`u $C>0$ est
une constante universelle.
\end{itemize}\end{lemme}
\textit{D\'emonstration.} En chaque point $t$ de $T$, d\'efinissons les $1$-jets 
$j_{0}(t)$, $\ldots,j_{n}(t)$ de section holomorphe de $T_t M$ dans $E$ par
\[ j_{0}(t)=J_1(s),\ j_{1}(t)=J_1(z_{1}s), \ldots ,\  j_{n}(t)=J_1(z_{n}s),\]
dans les coordonn\'ees $z$ centr\'ees en $t$ de \ref{voisinage bilipschitzien a une boule euclidienne}. 
Ces jets sont born\'es par $1$ et d'apr\`es la Proposition \ref{prolongement d'une famille de jets},
il existe des prolongements holomorphes $\sigma_{0},\ldots,\sigma_{n}:M\rightarrow E$
dont la norme est major\'ee par une constante $C$ ne d\'ependant que de $\delta$.

D'apr\`es le Lemme de Schwarz, nous avons $|\sigma_{0}|\geq 1/2$ sur la r\'eunion des boules $B_g(t,\varepsilon_1)$
centr\'ees en un point $t$ de $T$ et de rayon $\varepsilon_1>0$ ne d\'ependant que de $r$ et de $C$. 
En effet, le quotient $\sigma_{0}/s$ est une fonction holomorphe
\`a valeurs dans ${\bf C}$ d\'efinie sur la boule $B_g(t,r)$
et born\'ee par $C$.

Une autre version du Lemme de Schwarz montre que la fonction $f=(\frac{\sigma_{1}}{\sigma_{0}},\ldots,
\frac{\sigma_n}{\sigma_{0}})$ v\'erifie $||df_{x}||\geq 1/2$ pour tout point
$x$ de $B_g(t,\varepsilon_{2})$, o\`u $\varepsilon_{2}$ est un r\'eel ne d\'ependant que de $r$ et de $C$.
En effet, elle d\'efinit une application holomorphe
$B_g(t,r)\rightarrow {\bf C}^{n}$ born\'ee par $C$ et dont la diff\'erentielle
en $0$ est l'identit\'e. Le Lemme est d\'emontr\'e avec $\varepsilon=\min(\varepsilon _1, \varepsilon_2)$.

\begin{lemme} Soient $0<\varepsilon<\delta$ deux r\'eels. Il existe un nombre fini 
$T_{1},\ldots,T_{k}$ de parties $\delta$-s\'epar\'ees de $M$ dont la r\'eunion 
est $\epsilon$-dense.\end{lemme}
\textit{D\'emonstration.} Une partie $\varepsilon$-s\'epar\'ee et maximale 
pour l'inclusion est $\varepsilon$-dense. Prenons en une $T$. 
Choisissons une sous-partie $T_{1}$ de $T$
qui est $\delta$-s\'epar\'ee et maximale pour l'inclusion, puis une partie
$T_{2}\subset T-T_{1}$ $\delta$-s\'epar\'ee et maximale pour l'inclusion etc.
Ce proc\'ed\'e s'arr\^ete au bout d'un certain temps, c'est \`a dire
qu'il existe $l$ tel que $T=T_{1}\cup T_{2}\cup \ldots \cup T_{l}$.
En effet, supposons qu'il existe un point $t$ de $T$ qui ne soit dans aucun des
$T_{m}$, pour $m=1,\ldots,l$. Par maximalit\'e, il existe un point $t_{m}$ de $T_{m}$
dans $B_g(t,\delta)$ pour tout $m=1,\ldots,l$. Ces points sont deux \`a deux distincts
puisque les parties $T_{m}$ sont disjointes. Il y a donc $l+1$ points de $T$ dans
$B_g(t,\delta)$, et les $l+1$ boules de rayon $\varepsilon/2$ centr\'ees en ces points sont
deux \`a deux disjointes dans $B_g(t,\delta)$.
Nos hypoth\`eses de g\'eom\'etrie born\'ee montrent que $l+1$ est inf\'erieur \`a une 
constante ne d\'ependant que de la courbure de Ricci de $g$ et du rayon d'injectivit\'e. 
Le Lemme est d\'emontr\'e. 

\vspace{0.2cm}

Nous sommes en mesure d'achever la d\'emonstration du 
Th\'eor\`eme~\ref{immersion d'une variete non compacte}. 
Pour chacunes de ces parties $T_{m}$, 
consid\'erons les sections holomorphes $\sigma_{m,0},\ldots,\sigma_{m,n}$ de $E$ 
construites au Lemme~\ref{sections holomorphes pour l'immersion}
et l'application
\[\pi =[\sigma_{m,j}]_{1\leq m \leq l, 0\leq j\leq n}:M\rightarrow {\bf C}P^N,\]
o\`u $N=(n+1)l-1$. Compos\'ee avec les projections $p_{m}:{\bf C}P^N\rightarrow
{\bf C}P^n$ on retrouve les applications m\'eromorphes $p_{m}\circ \pi =\pi_{T_{m}}$.
On a donc l'in\'egalit\'e $\pi ^{*}\Omega_{FS}\geq C\Omega$ pour une constante strictement positive $C$. 

De plus en chaque point $x$ de $M$, l'une au moins des sections $\sigma_{m,j}$ est de 
norme sup\'erieure \`a $1/2$ et toutes sont uniform\'ement born\'ees par une constante universelle. 
La tir\'ee en arri\`ere par $\pi$ de la m\'etrique 
de Fubini-Study $\Omega_{FS}$ de ${\bf C}P^N$ s'exprimant par la formule 
\[ \pi^{*} \Omega_{FS}= \frac{\sqrt{-1}}{2} \overline{\partial}\partial
\log \big(  \sum_{l} |\sigma_{l}|^{2}\big),\]
nous obtenons l'in\'egalit\'e $\pi^{*}\Omega_{FS}\leq D\Omega$, pour une constante $D>0$.
Le Th\'eor\`e\-me~\ref{immersion d'une variete non compacte} est d\'emontr\'e.

\subsection{D\'emonstration du Th\'eor\`eme \ref{surfaces orientees}}

Soit $(\Sigma,g)$ une surface riemannienne orient\'ee 
\`a g\'eom\'etrie born\'ee est projective.
D'apr\`es le Th\'eo\-r\`eme d'Ahlfors-Bers \cite{A-B}, 
toute surface riemannienne orient\'ee est confor\-m\'ement plate 
et admet donc une structure de surface de Riemann. Dans 
le cas o\`u $\Sigma$ est compacte, le Th\'eor\`eme de Riemann montre que 
$\Sigma$ se plonge holomorphiquement dans ${\bf C}P(3)$. 
Si $\Sigma$ est non compacte, alors $H^2(\Sigma,{\bf R})$ est nul, 
et ceci implique l'existence d'un fibr\'e en 
droites holomorphe $E\rightarrow \Sigma$ muni d'une m\'etrique dont la courbure 
est la forme volume $v_g$. Cette m\'etrique est 
\`a g\'eom\'etrie born\'ee selon \ref{hypothese de geometrie bornee} si la courbure de $g$ 
est born\'ee uniform\'ement et le rayon d'injectivit\'e 
de $(\Sigma,g)$ est minor\'e uniform\'ement par une constante positive. 
Le r\'esultat d\'ecoule donc du Th\'eor\`eme~\ref{immersion d'une variete non compacte}. 

\subsection{Exemples} Un rev\^etement infini d'une vari\'et\'e projective hermitienne est un exemple 
de vari\'et\'e hermitienne qui s'immerge localement bilipschitziennement et holomorphiquement 
dans un espace projectif complexe. 

Notons que la famille des vari\'et\'es hermitiennes qui s'immergent localement bilipschitziennement et holomorphiquement 
dans un espace projectif complexe est stable par produit et par  
\'eclatement le long d'une partie s\'epar\'ee ou d'une sous-vari\'et\'e \`a g\'eom\'etrie born\'ee. 
Ceci donne donc d\'eja beaucoup d'exemples de vari\'et\'es hermitiennes projectives en dimension sup\'erieure.

Nous d\'ecrivons ci-apr\`es des exemples de vari\'et\'es quasi-projectives non compactes munies 
d'une m\'etrique hermitienne \`a g\'eom\'etrie born\'ee 
qui ne s'immergent pas bilipschitziennement et holomorphiquement dans un espace projectif complexe.

Rappelons que parmi les vari\'et\'es holomorphiquement parall\'elisables 
compactes du type $G/\Gamma$, o\`u $G$ est un groupe de Lie complexe et
$\Gamma$ un r\'eseau cocompact de $G$, celles qui sont k\"ahl\'eriennes sont 
les tores complexes ${\bf C}^n/\Gamma$. Soit $GA$ le groupe affine
complexe des matrices $2\times 2$ triangulaires sup\'erieures
de d\'eterminant $1$, muni d'une m\'etrique hermitienne invariante
\`a droite. 
\begin{proposition} Il n'existe pas d'immersion holomorphe
$\pi: GA \rightarrow M$ \`a valeurs dans une vari\'et\'e
k\"ahl\'erienne $(M,g)$ v\'erifiant l'in\'egalit\'e
\[
\frac{1}{C}\pi^{*}g \leq G \leq C
\pi^{*}g,\] o\`u $C$ est une constante
strictement positive et $G$ est une m\'etrique hermitienne sur $GA$ invariante \`a droite. 
\label{non immersion du groupe affine}\end{proposition} 
\textit{D\'emonstration.} Sur $GA$,
faisons agir le groupe \`a $1$ param\`etre $t$ complexe par
multiplication \`a gauche par les matrices
$$\left(\begin{array}{cc}
1 & t\\
0 & 1
\end{array}\right) .$$
Les orbites de ce groupe d\'efinissent un feuilletage ${\mathcal F}$
dont les feuilles sont paraboliques. Supposons qu'il existe une
immersion v\'erifiant les hypoth\`eses de la Proposition.

Commen\c{c}ons par construire un courant strictement positif
ferm\'e de type $(1,1)$ sur $M$ en utilisant une id\'ee d\^ue \`a L.~Ahlfors \cite{Ahlfors1}. 
Les feuilles de $\mathcal F$ sont des courbes enti\`eres sur lesquelles la restriction de $G$ 
induit la m\'etrique euclidienne. Prenons en une $F$ et consid\'erons les courants sur $M$
d\'efinis pour tout $r>0$ par
\[ T_r (w) = \frac{1}{\mathrm{Aire}(D_r)} \int_{D_r} \pi^{*}w,\]
o\`u $D_r$ est un disque euclidien de $F$ de rayon $r$, et l'aire \'etant mesur\'ee 
avec la m\'etrique $\pi^*g$. C'est une suite de courants
de norme $1$ sur $M$ et donc on peut en extraire une sous-suite
$T_{r_i}$ convergente vers un courant $T$. Comme l'immersion $\pi$ est localement 
bilipschitzienne la longueur du bord $\partial{D_r}$ est n\'egligeable devant
l'aire de $D_r$ pour la m\'etrique $\pi^*g$ lorsque $r$ tend vers l'infini, 
ce qui montre que $T$ est un courant ferm\'e.

Montrons que $T$ est homologue \`a $0$ dans $M$. Pour cela
remarquons que les transformations $\phi_{s}:GA\rightarrow GA$
obtenues par multiplication \`a gauche par les matrices
$$\left(\begin{array}{cc}
e^{-s} & 0\\
0 & e^{s}
\end{array}\right) $$
pr\'eservent le feuilletage ${\mathcal F}$ en contractant la restriction de $G$ aux feuilles 
par $e^{-2s}$. Nous avons donc des majorations
\begin{equation}\label{contraction} |\phi_{s}^{*} w |_{\pi ^* g,\infty}
\leq C e^{-2s}|w|_{\pi^*g,\infty} \end{equation} pour toute
$(1,1)$-forme $w$ d\'efinie le long des feuilles de ${\mathcal F}$. A
priori il n'y a pas de raison que $\phi_s$ s'\'etende \`a une
transformation de $M$. Cependant, imaginons un instant que ce soit
le cas. Soit $\omega$ la $(1,1)$-forme sur $M$ naturellement associ\'ee \`a $g$. 
Cette forme est ferm\'ee donc 
$T(\phi_s ^* \omega)=T(\omega)$ pour tout $s$, et
d'apr\`es \ref{contraction} nous avons $T(\phi_s ^* \omega)\leq
Ce^{-2s}$, ce qui montre que $T(\omega)=0$ et contredit la stricte
positivit\'e de $T$. En g\'en\'eral, $\phi_s$ ne se prolonge pas
\`a un flot sur $M$, mais si $w$ est l'image r\'eciproque par
$\pi^ *$ de $\omega$, nous prouvons que la limite
\[ \lim _i T_{r_i}(\phi_s ^* w)\]
existe et est \'egale \`a $T(\omega)$. Pour cela il suffit de
remarquer que $w$ est une $2$-forme ferm\'ee et born\'ee, et
diff\`ere de $\phi_{s}^{*} w$ de la diff\'erentielle d'une
$1$-forme \textit{born\'ee} $\eta$. En effet des relations de
Cartan nous d\'eduisons la formule
\[ \eta = d\big(\int_0 ^t i_X \phi_t^*w dt\big),\]
o\`u $X$ est le champ de norme $1$ induit par le flot $\phi_s$.
Nous avons alors
\[\lim_{i}\frac{1}{\mathrm{Aire}(\pi(D_{r_{i}}))}\int_{\partial{D_{r_{i}}}} \eta =0,\]
et nous en d\'eduisons la convergence de
$\frac{1}{\mathrm{Aire}(\pi(D_{r_{i}}))}\int_{D_{r_{i}}}
\phi_{s}^{*}w$ vers $T(\omega)$. D'apr\`es \ref{contraction} nous
avons
\[|T(\omega)|\leq C e^{-2s}.\]
Cela \'etant vrai pour tout $s$, c'est que $T(\omega)=0$. C'est
une contradiction puisque $T$ est strictement positif.

\vspace{0.2cm}

Notre m\'ethode montre aussi que la vari\'et\'e de Hopf
${\bf C} ^{2} -\{0\}$ muni de la m\'etrique $\frac{|dz|}{|z|}$
invariante par les homoth\'eties ne peut pas s'immerger localement
bilipshitziennement et holomorphiquement dans une vari\'et\'e
k\"ahl\'erienne. Elle peut s'adapter \`a bien d'autres espaces 
homog\`enes. Cependant nous ne savons pas si le groupe de Heisenberg 
complexe muni d'une m\'etrique invariante \`a droite 
s'immerge ou pas holomorphiquement et de mani\`ere bilipshitzienne dans 
un espace projectif complexe.
Nous reviendrons sur ce genre d'obstructions au paragraphe~\ref{cycle feuillete}.

\section{Continuit\'e des s\'eries fuchsiennes}\label{continuite des series fuchsiennes}

Dans cette partie nous \'etablissons des propri\'et\'es de continuit\'e des sections minimisantes, 
et de leurs s\'eries fuchsiennes associ\'ees.  
 
Soient $p_i: E_i\rightarrow M_i $, $i=1,2$ deux fibr\'es en droites holomorphes, et 
$|.|_i$ des m\'etriques hermitiennes lisses sur $E_i$ de courbure strictement positive 
et \`a g\'eom\'etrie born\'ee. Nous supposons que pour $i=1,2$, 
le rayon $r(|.|_i)$ est sup\'erieur \`a $r_0$ et 
que la courbure de Ricci de la m\'etrique k\"ahl\'erienne $g_i$ que $|.|_i$ induit 
sur $M_i$ est major\'ee par $1/4$. 
Sous ces conditions, tout $1$-jet $j:T_{t}M_i \rightarrow E_i$ de section holomorphe de $E_i$ se prolonge en 
une section holomorphe de $L^2(|.| _{i,t})$, o\`u 
$|.|_{i,t}= \exp(\alpha d(.,t)) |.|_i$ (voir Lemme \ref{section L2alphax}). 
Rappelons que $m(j)$ d\'esigne celle de norme minimale.

Soient $0<\varepsilon<1$ et $R>>1$ deux r\'eels.  
Supposons qu'il existe des domaines $D_i \subset M_i$ contenant des boules $B(x_i, R)$ 
et un isomorphisme de fibr\'e 
en droites complexes $\phi : p_1^{-1} (D_1) \rightarrow p_2^{-1} (D_2)$ tel que

\vspace{0.2cm}

(i) $\phi$ est $(1+\varepsilon)$-bilipschitzien,

\vspace{0.2cm}

(ii) $\exp(-\varepsilon) |.| \leq \phi _* |.| \leq \exp(\varepsilon) |.|$,

\vspace{0.2cm}

(iii) Pour tout $1$-jet $\omega_1$ (resp. $\omega_2$) de section holomorphe de $E_1$ (resp. $E_2$), 
$|\overline{\partial} \phi_* \omega_1 | \leq \varepsilon |\omega_1|$ 
(resp. $|\overline{\partial} \phi^* \omega_2 |\leq \varepsilon |\omega_2|$).    

\vspace{0.2cm}

\begin{lemme}[Continuit\'e des sections minimisantes] 
Soient $t_1$ un point de $B(x_1,R/8)$ et $t_2 = \phi(t_1)$. 
Supposons que $\phi$ est holomorphe sur la boule $B(x_1,1)$. 
Soit $j_1$ un $1$-jet de section holomorphe de $E_1$ en $t_1$ et $j_2 = \phi_* j_1$. Alors 
\[  |\phi m(j_1)(x_1) - m(j_2)(x_2) |\leq C\exp(-\alpha d(x_1,t_1) ) f(R,\varepsilon) |j_1|,\]
o\`u $f(\varepsilon,R)$ est une fonction qui tend vers $1/\sqrt{R}$ lorsque $\varepsilon$ tend vers $0$,
et $C$ est une constante universelle. 
\label{lemme de Gromov}\end{lemme}

Les sections minimisantes passant par un jet en un point $t$ 
ne d\'epen\-dent donc essentiellement que de la g\'eom\'etrie de $|.|$ en restriction 
\`a une grande boule centr\'ee en $t$. Cette propri\'et\'e a \'et\'e mise en \'evidence par M.~Gromov
(\cite{Gromov1}, p. 387, Remark (a) of 3.3.5).

\vspace{0.2cm}

\textit{D\'emonstration.}  
Observons que puisque $\varepsilon <1$, $\phi $ est $2$-bilipschitzien. Nous avons donc les inclusions~:
\[ B(t_2, R/2)\subset B(x_2, 3R/4) \subset B(x_2, R).\] 
Consid\'erons une fonction lisse $\psi: M_2 \rightarrow [0,1]$ qui vaut $1$ sur $B(t_2,R/2)$,
$0$ \`a l'ext\'erieur de $B(x_2,R)$, et telle que \[ |d\psi |\leq C/R,\] 
la constante $C$ ne d\'ependant que de $r_0$. 
Nous voulons perturber la section \[\psi \phi_* m(j_1):M_2\rightarrow E_2\] 
en une section holomorphe $h_2:M_2\rightarrow E_2$ dans la norme $|.|_{t_2,2}$.\footnote{
Pour tout $U$ dans $\mathcal U$, toute partie $B\subset M_i$, tout point $z$ de $M$, 
et toute section $\tau :B\rightarrow E$, notons 
\[  |\tau|_{z,B,2}= \sqrt{\int_B  |\tau(y)|_z^2 dv_g(y)}\ \ \ \mathrm{et}\ \ \ \  |\tau|_{z,2} = |\tau |_{z,M,2}.\]}

\vspace{0.2cm} 

\textit{Comportement des m\'etriques $|.|_{t}$.}
Les propri\'et\'es suivantes portant sur les m\'etriques 
$ |.|_{i,t_i}$ d\'ecoulent directement des propri\'et\'es (ii) et (iii)~: 
il existe une fonction $\eta(R,\varepsilon)$  
qui tend vers $0$ lorsque $\varepsilon$ tend vers $0$ et telle que 

\vspace{0.2cm}

(ii)' $\exp(-\eta ) |.|_{t_2}\leq \phi_* |.|_{t_1}\leq \exp(\eta) |.|{t_2}$,

\vspace{0.2cm} 

(iii)' Pour tout $1$-jet $\omega_1$ (resp. $\omega_2$) de section holomorphe de $E_1$ (resp. $E_2$), 
$|\overline{\partial} \phi_* \omega_1 |_{t_2} \leq \eta |\omega_1|_{t_1}$ 
(resp. $|\overline{\partial} \phi^* \omega_2 |_{t_1}\leq \eta |\omega_2|_{t_2}$).

\vspace{0.4cm}

\textit{Estimation du $\overline{\partial}$.} Nous avons 
\[ | \overline{\partial }(\psi \phi m(j_1))|_{t_2,2}\leq | (\overline{\partial }\psi) \phi m(j_1))|_{t_2,2}+
| (\overline{\partial }\phi m(j_1))\psi|_{t_2,2}.\]
Commen\c{c}ons par majorer le premier terme. Nous avons 
\[  |(\overline{\partial} \psi) \phi m(j_1)|_{t_2,2}^2 \leq \big(\frac{C}{R}\big)^2 |\phi m(j_1)|_{t_2,B(x_2,R)}^2\]
\[ = \big(\frac{C}{R}\big)^2 \int_{B(x_2,R)} |\phi m(j_1)|_{t_2}^2 dv_{g_2}\leq 
\big(\frac{C}{R}\big)^2 |m(j_1)|_{t_1,2}^2.\]
D'apr\`es le Lemme~\ref{section L2alphax}, nous obtenons donc
\begin{equation}\label{premiere equation}  
|(\overline{\partial} \psi) \phi m(j_1)|_{t_2,2} \leq \frac{C}{R}|j_1|.\end{equation} 
Majorons le deuxi\`eme terme. D'apr\`es (iii)', nous avons 
\[ |\psi \overline{\partial } \phi m(j_1)|_{t_2,2} ^2 \leq |\overline{\partial} \phi m(j_1)|_{2_1,2}
\leq \int_{B(x_2,R)} |\overline{\partial} \phi m(j_1)|_{t_2} ^2 dv_{g_2} \]
\[\leq C\eta ^2  |J_1 m(j_1)|_{t_1}^2,\]
o\`u $J_1 m(j_1)$ d\'esigne le premier jet de $m(j_1)$. Les in\'egalit\'es de Garding 
nous donnent en tout point $y_1$ de $M_1$,
\[  |J_1 m(j_1)(y_1)|_{t_1} ^2 \leq C \int _{B(y_1,1)} |m(j_1)(z_1)|_{t_1}^2dv_{g_1}(z_1),\]
donc 
\[ |J_1 m(j_1)|_{2}^2\leq C\int _{d(y_1,z_1)\leq 1} |m(j_1)(z_1)|_{t_1} ^2 dv_{g_1}(y_1)dv_{g_1}(z_1)\]
\[ = C\int_{z_1\in M_1} \mathrm{vol}(B(z_1,1)) |m(j_1)(z_1)|_{t_1}^2 dv_{g_1}(z_1) \leq C|m(j_1)|_{t_1}^2.\]
Nous en d\'eduisons 
\[  |\psi \overline{\partial} \phi m(j_1)|_{t_2,2}\leq C\eta |j_1|,\]
et en combinant avec \ref{premiere equation}, 
\[ |\overline{\partial} (\psi \phi m(j_1))|_{t_2,2} \leq C(\eta + \frac{1}{R}) |j_1|.\]
\vspace{0.1cm}

\textit{Perturbation.} 
Le Lemme~\ref{section L2alphax} fournit une section holomorphe $h_2:M_2\rightarrow E_2$ telle que 
\begin{equation}\label{seconde equation} |h_2-\psi \phi m(j_1)|_{t_2,2} 
\leq C(\eta + \frac{1}{R}) |j_1|.\end{equation}
L'in\'egalit\'e de Garding 
donne, puisque $h_2-\psi \phi m(j_1)$ est holomorphe sur $B(x_2,1/2)$~: 
\[ |J_1 h_2(t_2) -j_2|  \leq C(\eta + \frac{1}{R})|j_1|.\]
Quitte \`a ajouter \`a $h_2$ la section $m(j_2 -J_1 h_2))$, nous pouvons supposer que l'on a en plus 
\[  J_1 h_2 = j_2.\]
\vspace{0.1cm}

\textit{Pythagore et une in\'egalit\'e.} Comme $m(j_2)$ est la section holomorphe de 
$L^2(|.|_{t_2})$ passant par $j_2$ de norme minimale, 
elle est orthogonale \`a l'espace des sections dont le premier jet en $t_2$ s'annule. Nous avons donc
\[|h_2|_{t_2,2}^2 = |m(j_2)|_{t_2,2}^2+|h_2-m(j_2)|_{t_2,2}^2\geq |m(j_2)| _{t_2,2}^2.\]
En vertu de~\ref{seconde equation}, on trouve 
\[ |\psi \phi m(j_1)|_{t_2,2} \geq ( 1-C(\eta + \frac{1}{R})) |m(j_2)|_{t_2,2}.\]
Or 
\[  |\psi \phi m(j_1)|_{t_2,2} \leq |\phi m(j_1)|_{t_2,D_2,2}\leq (1+\varepsilon)^{n+3/2}|m(j_1)|_{t_1,D_1,2},\]
d'o\`u les in\'egalit\'es
\begin{equation} \label{inegalite} (1+\varepsilon)^{n+3/2} |m(j_1)|_{t_1,2} \geq |\psi\phi m(j_1)|_{t_2,2}\geq 
( 1-C(\eta + \frac{1}{R}))  |m(j_2)|_{t_2,2}.\end{equation}  
\vspace{0.1cm}

\textit{Sym\'etrie.} Les hypoth\`eses du Lemme sont sym\'etriques. Les 
in\'egalit\'es \ref{inegalite} sont donc valables 
pour $m(j_2)$ aussi. Nous obtenons alors 
\[  (1+\varepsilon )^{n+3/2}|m(j_2)|_{t_2,2} \geq ( 1-C(\eta + \frac{1}{R}))  |m(j_1)|_{t_1,2}.\]
En r\'eutilisant \ref{inegalite} et \ref{seconde equation}, nous trouvons une fonction $\eta'(R,\varepsilon) $ 
qui tend vers $0$ lorsque $\varepsilon$ tend vers $0$ et telle que 
\[ ( 1+C(\eta' + \frac{1}{R})) |m(j_2)|_{t_2,2}\geq |h_2|_{t_2,2} \geq 
( 1-C(\eta' + \frac{1}{R})) |m(j_2)|_{t_2,2}.\]
Finalement, en vertu de l'\'egalit\'e de Pythagore,
\[ |h_2-m(j_2)|_{t_2,2} = \sqrt{|h_2|_{t_2,2} ^2-|m(j_2)|_{t_2,2}^2}\leq C\sqrt{ \eta '+ \frac{1}{R}}|j_2|\]
et de nouveau \`a cause de \ref{seconde equation}, on obtient
\[ |\phi m(j_1)-m(j_2)|_{t_2,B(t_2,R/2),2}\leq |\psi \phi m(j_1) -m(j_2)|_{t_2,2}\]
\[ \leq C\sqrt{\eta' + \frac{1}{R}}|j_2|.\]
Le Lemme~\ref{lemme de Gromov} d\'ecoule alors de l'in\'egalit\'e de Garding, en choisissant la fonction 
\[f(R,\varepsilon ) = C (\eta'(R,\varepsilon) + \frac{1}{R})^{1/2}.\]  

\vspace{0.2cm}

Nous nous int\'eressons \`a pr\'esent \`a la continuit\'e des s\'eries fuchsiennes. 
Pour assurer leur convergence, nous supposons que la courbure de Ricci des m\'etriques $g_i$ est major\'ee par le 
r\'eel $c_0$ donn\'e par le Lemme~\ref{convergence des series fuchsiennes}.
Rappelons que le r\'eel $\delta _0>>1$ 
a \'et\'e d\'efini \`a la Proposition~\ref{prolongement d'une famille de jets}. 

\begin{lemme} [Continuit\'e des s\'eries fuchsiennes] 
Soient $T_i$ des parties $\delta_0$-s\'epar\'ees de $M_i$,
et $j_i = \{ j_i(t)\} _{t\in M_i}$ des familles born\'ees de $1$-jets de 
section holomorphe de $E_i$ dont le support est contenu dans $T_i$. 
Supposons que $\phi$ est holomorphe sur chaque boule $B(t,1)$, o\`u $t\in T_1 \cap D_1$ et 
que $\phi
(T_1 \cap D_1) = T_2 \cap D_2$. Alors 
\[  |\phi \sigma (j_1) (x_1) -\sigma (j_2) (x_2) | \leq  
C (f(R,\varepsilon) \max (|j_1|_{\infty},|j_2|_{\infty}) +|j_2 - \phi_*j_1|_{\infty,D_2} ),\]
o\`u $f(R,\varepsilon)$ est une fonction qui tend vers $1/\sqrt{R}$ lorsque $\varepsilon$ tend vers $0$
et $C$ est une constante universelle. 
\label{lemme de continuite sur les series fuchsiennes}\end{lemme} 

\textit{D\'emonstration.} 
Nous avons 
\[ \phi(\sigma(j_1)(x_1)) -\sigma(j_2)(x_2) = \]
\[\sum_{t\in T_1\cap B(x_1,R/8)} \phi(m(j_1(t))(x_1))-m(j_2(\phi(t)))(x_2)\]
\[+\sum_{t\in M_1-B(x_1,R/8)} \phi(m(j_1(t))(x_1))-\sum_{t\in M_2-\phi(B(x_1,R/8)} m(j_2(t))(x_2).\]
Nous allons majorer la norme de chacun des termes de droite de cette \'egalit\'e. 

\vspace{0.2cm}

\textit{Majoration de la norme du premier terme.} Nous avons 
pour tout $t\in T_1 \cap B(x_1, R/8)$, 
\[  |\phi (m(j_1(t))(x_1)) - m(j_2(\phi(t)))(x_2)|\leq \]
\[| \phi(m(j_1(t))(x_1)) - m(\phi_* j_1(t))(x_2)| + | m(\phi_* j_1(t)-j_2(\phi(t)))(x_2)| .\] 
Ainsi, en appliquant le Lem\-me~\ref{lemme de Gromov},
\[ |\sum_{t\in T_1 \cap B(x_1,R/8)}\phi(m(j_1(t))(x_1))-m(j_2(\phi(t)))(x_2)| \leq \]
\[ f(R,\varepsilon) |j_1|_{\infty}  \sum_{t\in T_1 \cap B(x_1,R/8)} e^{-\alpha d(x_1,t)}+\]
\[ C|\phi_*j_1-j_2|_{\infty, D_1} \sum _{t\in T_2 \cap D_2} e^{-\alpha d(x_2,t)}.\]
Comme $T_1$ et $T_2$ sont des parties $\delta_0$-s\'epar\'ees, 
les techniques de~\ref{convergence des series fuchsiennes} 
montrent que les deux sommes apparaissant dans le terme de droite
de l'in\'egalit\'e pr\'ec\'edente sont major\'ees par 
une constante ne d\'ependant que de $r_0,\ \alpha,\ \delta_0$ et $c_0$. Nous obtenons donc 
l'in\'egalit\'e voulue
\[ |\sum_{t\in T_1\cap B(x_1,R/8)}\phi(m(j(t))(x_1))-m(j(\phi(t)))(x_2)| \leq \]
\[ C(f(R,\varepsilon) \max(|j_1|_{\infty} ,|j_2|_{\infty})+ |\phi_*j_1-j_2|_{\infty, D_1}).\]  
\vspace{0.1cm}

\textit{Majoration de la norme des deux autres termes.} Puisque les sections minimisantes 
sont \`a d\'ecroissance exponentielle pour la norme $|.|$,
nous avons 
\[ |\sum_{t\in T_1-B(x_1,R/8 )} \phi(m(j_1(t))(x_1))| \leq C|j_1|_{\infty} 
\sum_{t\in T_1-B(x_1,R/8)} e^{-\alpha d(x_1,t)}.\]
Les techniques de \ref{prolongement d'une famille de jets} montrent que 
\[ \sum_{t\in T_1-B(x_1,R/8)} e^{-\alpha d(x_1,t)} 
\leq C \int_{R/8 - \delta _0} ^{\infty} e^{-\alpha r/2} dr
\leq Ce^{-\frac{\alpha}{16/R}}\leq C/\sqrt{R},\]
la constante ne d\'ependant que de $r_0,\ \alpha,\ c_0$, et $\delta_0$. Nous obtenons donc 
\[ |\sum_{t\in T_1- B(x_1,R/8)} \phi(m(j_1(t))(x_1))| \leq C|j_1|_{\infty}/\sqrt{R}.\]
Pour le troisi\`eme terme, le m\^eme raisonnement donne 
\[  |\sum_{t\in M_2 -\phi(B(x_1,1/8\varepsilon))} m(j_2(t))(x_2)| \leq C|j_2|_{\infty}/\sqrt{R},\]
puisque 
\[M_2-\phi(B(x_1,R/8)) \subset M_2-B(x_2,R/16).\]
Le Lemme~\ref{lemme de continuite sur les series fuchsiennes} est d\'emontr\'e.

\section{S\'eries fuchsiennes sur une lamination}
\label{series fuchsiennes sur une lamination}
Une \textit{lamination par vari\'et\'es complexes} 
d'un espace topologique $X$ est un atlas ${\mathcal L}$ d'hom\'eomorphismes
$\varphi : U\rightarrow B\times T$ \`a valeurs
dans le produit d'une boule $B$ de ${\bf C}^{n}$ par 
un espace topologique $T$, appel\'e \textit{espace transverse},
tels que les changements de cartes soient de la forme
\[ (z,t)\mapsto (z'(z,t),t'(t)),\]
o\`u $z'$ d\'epend holomorphiquement de $z$. Les feuilles de $\mathcal L$ localement d\'efinies par $B\times t$ 
sont des vari\'et\'es complexes et forment une r\'eunion disjointe de $X$.

Sur une lamination par vari\'et\'es complexes $\mathcal L$, 
soit $\mathcal O$ le pr\'efaisceau des fonctions 
\textit{continues} $f:U\rightarrow {\bf C}$ d\'efinies sur un ouvert $U$ de $X$, 
qui sont holomorphes le long des feuilles de $\mathcal L$. 
Ce faisceau nous permet d'\'etendre aux laminations par 
vari\'et\'es complexes les notions analytiques classiques dont on dispose 
sur les vari\'et\'es complexes. Par exemple un fibr\'e en droites
holomorphe sur $\mathcal L$ est un \'el\'ement de \\$H^1(X,\mathcal O^*)$.

Consid\'erons aussi le pr\'efaisceau $C^{\infty}(\mathcal L)$ des fonctions $\varphi :U\rightarrow {\bf R}$  
d\'efinies sur un ouvert $U$ de $X$ qui sont lisses le long des feuilles, 
et telles que les d\'eriv\'ees \`a tout ordre par rapport 
aux coordonn\'ees $z$ soient aussi continues. Comme les changements de coordonn\'ees sont holomorphes, 
cette notion ne d\'epend pas du syst\`eme de coordonn\'ees choisies. 
Comme dans le cas des vari\'et\'es, il est possible de construire des partitions de l'unit\'e
associ\'ees \`a un recouvrement localement fini par des ouverts de $X$ avec des fonctions de $C^{\infty}(\mathcal L)$. 
Ceci permet de construire des m\'etriques lisses sur un fibr\'e en droites holomorphe au dessus de $\mathcal L$.

\subsection{Structure produit}
En g\'en\'eral, si $F$ est une feuille d'une lamination lisse, 
il n'est pas possible de d\'eformer un domaine~\footnote{Un domaine est un 
ouvert connexe relativement compact \`a bord lisse.} $D\subset F$ 
dans une feuille voisine. L'obstruction est le fait que la repr\'esentation d'holonomie 
\[  Hol : \pi_1(D,*) \rightarrow \mathrm{Homeo}( T,*) \] 
n'est pas triviale, o\`u $T$ est un espace transverse au point $*\in D$ 
et $\mathrm{Homeo}(T,*)$ d\'esigne le groupe des germes d'hom\'eo\-morphismes de $T$ qui fixent $*$.

\begin{lemme} 
Soit $D\subset F$ un domaine contenu dans une feuille $F$ pour lequel la repr\'esentation d'holonomie 
$Hol $ est triviale. Alors il existe un voisinage $\mathcal V$ de $D$ dans $X$ et un diff\'eomorphisme 
de lamination $\phi: \overline{D}\times T\rightarrow \mathcal V$. 

De plus, supposons que $\{ \overline{B_i} \}_{i\in I}$ soit une famille finie de boules ferm\'ees lisses 
disjointes contenues dans $D$, et que 
\[  \phi _i : \overline{B_i}\times T_i \rightarrow \mathcal V_i\]
soient des cartes feuillet\'ees avec $\phi_i (.,t_i)= id|_{\overline{B_i}}$.
Alors quitte \`a restreindre les $T_i$ il existe des identifications naturelles $T_i = T_j$ donn\'ees par l'holonomie, 
et l'on peut choisir le diff\'eomorphisme $\phi: D\times T\rightarrow \mathcal V$ en sorte que 
\[ \phi |_{\overline{B_i}\times T_i} = \phi_i.\]  
\label{deformation}\end{lemme} 

\textit{D\'emonstration.} Consid\'erons un plongement lisse $\pi : D'\rightarrow {\bf R}^N$ d'un voisinage 
$D'$ de $\overline {D}$ pour lequel la repr\'esentation d'holonomie est triviale. 
En utilisant une partition de l'unit\'e, il est possible de prolonger $\pi $ 
en une application d\'efinie sur un voisinage 
$V$ de $\overline{D}$ dans $X$. 
Nous la notons encore $\pi : V\rightarrow {\bf R}^N$ pour simplifier les notations. 

Soit $D'$ un domaine contenu dans $F$ 
tel que $\overline{D} \subset D' \subset \overline{D'}\subset V$. 
D'apr\`es le Th\'eor\`eme du voisinage tubulaire, il existe un voisinage $W$ de $\pi (\overline{D'})$ 
et une fibration lisse
\[  B \rightarrow W \stackrel {p}{\rightarrow } \overline{D'},\]
telle que $p\circ \pi | _{\overline{D'}}$ est l'identit\'e. 
Quitte \`a restreindre l'ouvert $V$, nous pouvons supposer que $p\circ \pi$ est une submersion 
le long des feuilles contenues dans $V$. Posons alors $U= V\cap (p\circ \pi)^{-1}(\overline{D'})$. 

En restriction \`a une feuille $F'$ 
de $U$, l'application $p\circ \pi:F' \rightarrow \overline{D'}$ est un rev\^etement 
dont le groupe est donn\'e par la repr\'esentation d'holonomie. Puisque cette derni\`ere est triviale,
l'application $p\circ \pi_{F}  $ est injective en restriction \`a une feuille de $U$. 

D'autre part, la restriction $p\circ \pi$ est surjective en restriction \`a des feuilles $F'$
de $U$ proches de $\overline{D'}$ puisque l'holonomie est d\'efinie sur tout domaine compact contenu 
dans une feuille. La premi\`ere partie du Lemme est d\'emontr\'ee.

La deuxi\`eme partie du Lemme r\'esulte du fait suivant. Soit $M$ une vari\'et\'e \`a bord et 
$B_i\subset M$ une famille finie de boules \`a bord lisse ne rencontrant pas le bord de $M$. 
Alors, si $\phi_i:\overline{B_i} \rightarrow M$ sont des plongements lisses assez proches de 
l'injection naturelle dans la topologie $C^1$, il existe un diff\'eomorphisme $\phi:M\rightarrow M $
assez proche de l'identit\'e dans la topologie $C^1$ tel que $\phi|_{\overline{B_i}} =\phi_i$.

\begin{lemme}
Soit $p: E \rightarrow \overline{D} \times T$ un fibr\'e en droites complexes lisse et $t_0 \in T$. 
Il existe un voisinage $T'$ de $t_0$ dans $T$ et un isomorphisme de fibr\'e en droites lisses
\[  \phi : E|_{\overline{D}\times \{t_0\}} \times T' \rightarrow E\]
qui induit l'identit\'e sur la base. 

De plus, avec les notations du Lemme \ref{deformation}, si l'on a des isomorphismes de fibr\'es 
en droites lisses 
\[  \phi_i : E|_{B_i \times \{t_0\}} \times T \rightarrow E|_{B_i\times T},\]
induisant l'identit\'e sur la base, il est possible de choisir $\phi$ en sorte que 
\[  \phi |_{B_i \times \{t_0\} } \times T' = \phi_i|_{B_i \times T'} .\]

\label{fibre en droites sur la lamination triviale}\end{lemme}

\textit{D\'emonstration.} 
Il existe un voisinage $T"$ de $t_0$ dans $T$, un recouvrement $\{U_j\}$ de $D\times T"$ 
par des ouverts $U_j$ de la forme $B_j \times T"$, pour lesquels il existe des 
isomorphismes de fibr\'es en droites lisses 
\[  \phi_j : E|_{B_j \times \{t_0\} } \times T" \rightarrow E|_{B_j \times T"},\]
qui induisent l'identit\'e sur la base et tels que 
\[  \phi_j |_{E|_{B_j \times \{t_0\}}\times \{t_0\}} = id.\]
Soit $\{ \chi _j\}_j$ une partition de l'unit\'e associ\'ee au recouvrement $\{U_j\}_j$.
\'Etendons $\chi_j \phi_j$ \`a un endomorphisme de fibr\'e en droites lisses
\[  \widetilde{\chi_j \phi_j} : E|_{\overline{D}\times \{t_0\}} \times T" \rightarrow E|_{\overline{D}\times T"},\]
et induisant l'identit\'e sur la base. La somme 
\[  \phi = \sum_j  \widetilde{\chi_j \phi_j}\]
est un endomorphisme de fibr\'e en droites lisse qui vaut l'identit\'e au dessus de $\overline{D}\times \{t_0\}$, 
et qui est l'identit\'e sur la base~: il existe donc un voisinage $T'$ de $t_0$ dans $T"$ tel que
\[  \phi : E|_{\overline{D}\times \{t_0\}} \times T' \rightarrow E,\]
est un isomorphisme de fibr\'e en droite lisse induisant l'identit\'e sur la base. 

Pour la deuxi\`eme partie, on commence par \'etendre les isomorphismes 
\[  
\phi_i : E|_{B_i \times \{t_0\}} \times T \rightarrow E|_{B_i\times T},\]
\`a des isomorphismes de fibr\'es en droites lisses 
\[  
\phi_i : E|_{B'_i \times \{t_0\}} \times T \rightarrow E|_{B'_i\times T},\]
induisant l'identit\'e sur la base, o\`u $B'_i$ est une boule ferm\'ee contenant 
$B_i$ dans son int\'erieur. On compl\`ete la famille $B'_i \times T"$ par des ouverts 
$U_j$ comme pr\'ec\'edemment qui ne rencontrent pas $B_i \times T$. Il suffit alors de 
prendre une partition de l'unit\'e telle que $\chi_i$ est identiquement constante \'egale \`a $1$ 
sur $B_i \times T$. Le lemme est d\'emontr\'e.

\subsection{Tubes }

Un \textit{bout de transversale} est une partie ferm\'ee $\tau$ de $X$ pour laquelle il existe une famille 
de bo\^{\i}tes $B\times T$ recouvrant $X$ telles que $\tau \cap (B\times T) \subset \{0\} \times T$.
Localement, un bout de transversale se plonge dans l'espace transverse $T$. On peut donc d\'efinir son bord 
pour la topologie transverse. 

Consi\-d\'e\-rons un bout de transversale $\tau$, et la r\'eunion disjointe 
$X_{\tau}$ des rev\^etements universels $\widetilde{F_t}$ des feuilles passant par un point $t$ de $\tau$. 
Cet ensemble est muni d'une topologie d\'efinie de la fa\c{c}on suivante~: un ouvert est l'ensemble des classes 
d'homotopies de lacets qui peuvent \^etre repr\'esent\'ees par un chemin contenu dans un ouvert donn\'e 
de $X$. L'espace topologique $X_{\tau} $ est appel\'e le \textit{tube} de $\tau$. 
Il est bien connu que l'espace $X_{\tau}$ est Hausdorff si et seulement si 
les feuilles passant par $ \tau$ ne supportent pas de \textit{cycle \'evanouissant} (voir par exemple \cite{Brunella,Moore-Schochet}). 

\begin{definition} Un \textit{cycle \'evanouissant} de $\mathcal L$ 
est un lacet $\gamma$ contenu dans une feuille et non homotope \`a 
un point dans sa feuille, mais qui est limite de lacets $\gamma _{n}$ contenus dans des feuilles
voisines et qui, eux, sont homotopes \`a un point dans leur feuille. 
\end{definition}

L'exemple typique d'un cycle \'evanouissant est un cylindre bouch\'e d'un cot\'e par un disque qui 
accumule \`a l'autre extr\'emit\'e sur un tore.

Lorsque $X_{\tau}$ est Hausdorff, 
il h\'erite de la structure d'une lamination par vari\'et\'es complexes, induite par la projection naturelle 
$r: X_{\tau} \rightarrow X$ qui est un hom\'eo\-morphisme local. La projection 
$X_{\tau}\rightarrow X$ est alors un biholomorphisme local de lamination par vari\'et\'es complexes.

L'avantage des tubes est le fait que pour tout domaine contenu dans l'une de leur feuille 
la repr\'esentation d'holonomie est triviale. On peut donc leur appliquer le Lemme \ref{deformation}.

\subsection{Construction de s\'eries fuchsiennes continues}

Soit $\mathcal L$ une lamination par vari\'et\'es complexes d'un espace compact $X$, 
$E\rightarrow \mathcal L$ un fibr\'e en droites holomorphe muni d'une m\'etrique hermitienne $|.|$ 
lisse dont la courbure en restriction \`a chaque feuille est strictement positive. 
Notons $g$ la m\'etrique k\"ahl\'erienne lisse le long des feuilles induites par la courbure 
de $|.|$.
Ces m\'etriques sont automatiquement \`a g\'eom\'etrie born\'ee parce que l'espace total est compact. 
Plus pr\'ecis\'ement il existe des r\'eels $r>0$ et $c>0$ tels 
que $r(|.|_{|F})\geq r$ et la courbure de Ricci de $g_F$
est born\'ee par $c$ pour toute feuille $F$. Nous supposons que $r\geq r_0$ et $c\leq c_0$, o\`u 
le r\'eel $r_0$ a \'et\'e d\'efini au Lemme~\ref{section L2alphax} 
et le r\'eel $c_0$ au Lemme~\ref{convergence des series fuchsiennes}. Ce n'est pas une grande restriction 
quitte \`a consid\'erer des puissances assez grandes de $E$.

Soit $j=\{ j_x\}_{x\in X}$ une famille \textit{born\'ee} de $1$-jets 
de sections holomorphes de $E$, dont le support est un bout de transversale $\mathcal T$ induisant 
sur chaque feuille une partie $\delta_0$-s\'epar\'ee. 
Sur le rev\^etement universel $r_F : \widetilde {F} \rightarrow F$ d'une feuille $F$ de $\mathcal L$, 
la s\'erie fuchsienne $\sigma ( r_F ^* j_F)$ converge et est invariante par le groupe du rev\^etement. 
Elle d\'efinit donc une section holomorphe de $E$ en restriction \`a $F$. La collection de toutes ses sections 
est une section born\'ee $\tilde{\sigma}(j):X\rightarrow E$ holomorphe sur chaque feuille. 
C'est la \textit{s\'erie fuchsienne} associ\'ee \`a la famille $j$. La raison pour laquelle on la note diff\'eremment 
est le fait que ces s\'eries ne sont pas d\'efinies directement sur les feuilles mais sur leur rev\^etement 
universel.

\begin{lemme} Supposons que $\mathcal L$ n'a pas de cycle \'evanouissant et que 
la famille $j$ est continue et s'annule sur le bord de $\mathcal T$. Alors $\tilde{\sigma}(j)$ 
est continue.\label{continuite des series fuchsiennes laminees}\end{lemme} 

\textit{D\'emonstration.} 
Soit $x$ un point de $X$ o\`u l'on veut d\'emontrer la continuit\'e de $\tilde{\sigma}(j)$. 
Consid\'erons un bout de transversale $\tau $ passant par $x$ et le tube $X_{\tau}$~: 
c'est une lamination par vari\'et\'es complexes puisque $\mathcal L$ n'a pas de cycle \'evanouissant. 
Enfin, donnons nous une section locale ne s'annulant pas $s: U\rightarrow r^*E$ d\'efinie dans un voisinage 
de $x$ et \'ecrivons $\tilde{\sigma}(j)= fs$.

Soit $R>>1$ et $D\subset \widetilde{F_x}$ un domaine contenant la boule 
$B(x,R)$. L'image r\'eciproque $r^{-1} (\mathcal T)$
est un bout de transversale de $X_{\tau}$ qui intersecte 
$D$ sur une partie $\delta_0$-s\'epar\'ee. Pour tout $t \in r^{-1}(\mathcal T)\cap D$, 
soient $\phi_t :B \times T_t\rightarrow \mathcal V_t$
des cartes locales holomorphes, o\`u $B$ est la boule unit\'e de ${\bf C }^n$, $T_t$ est un espace transverse 
contenant $t$, et $\mathcal V_t$ un voisinage de $t$. Soient aussi 
$s_t: \mathcal V_t \rightarrow r^* E$ des sections locales holomorphes trivialisantes. 
On peut supposer que l'image $\phi_t (B\times t)$ contient la boule $B(t,1)$. 
   
Appliquons les Lemmes \ref{lemme de continuite sur les series fuchsiennes}, \ref{deformation} et 
\ref{fibre en droites sur la lamination triviale}~: nous en d\'eduisons qu'il existe un voisinage 
$x\in V\subset U$ de $x$ et une structure de bo\^{\i}te lisse $V=B\times T$ (provenant de 
la structure produit $D\times T$ d'un voisinage de $D$), avec $x=(0,t_0)$ telle que 
\[  \liminf _{t\rightarrow t_0} |f(.,t)-f(.t_0)| \leq C\frac{|j|_{\infty}}{\sqrt{R}}.\]
Les in\'egalit\'es de Garding montrent que $f$ est uniform\'ement 
continue le long des feuilles~; cela prouve que $f$ est continue, et \`a fortiori $\tilde{\sigma}(j)$.

\begin{remarque} 
Lorsque la lamination est transversalement lipshitzienne, que le fibr\'e en
droites $E$ est lipshitzien, et que $j$ est lipschitzienne, 
le module de continuit\'e des s\'eries fuchsiennes $\tilde{\sigma}(j)$ n'est \`a priori major\'e que par 
\[ \delta(\sigma(j))(\varepsilon) \leq C\sqrt{\frac{1}{\log 1/\epsilon}}.\]
Nous n'obtenons donc que tr\`es peu de r\'egularit\'e
transverse avec cette m\'ethode.\end{remarque}

\section{Immersion de laminations par vari\'et\'es complexes}\label{immersion projective d'une lamination}

\subsection{D\'emonstration du Th\'eor\`eme \ref{theoreme de Kodaira lamine}}\label{demonstration du theoreme de Kodaira lamine}
Soit $\mathcal L$ une lamination par vari\'et\'es complexes 
d'un espace compact $X$ qui n'a pas de cycle \'evanouis\-sant, et $E\rightarrow \mathcal L$ un fibr\'e en droites holomorphe 
muni d'une m\'etrique \`a courbure strictement positive le long des feuilles. Il s'agit de montrer que 
$\mathcal L$ est projective.

Nous construisons des applications  
\[ \pi : X\rightarrow {\bf C}P^N \]
qui immergent holomorphiquement les feuilles de $\mathcal L$, et qui s\'eparent deux points distincts donn\'es
$x_1$ et $x_2$ de $X$. 
Les applications $\pi$ sont de la
forme $[\sigma_{0},\ldots,\sigma_{N}]$, o\`u les $\sigma_{l}$ sont
des sections d'une puissance assez grande $E^{\otimes k}$ de $E$ donn\'ees par le lemme suivant. 

\begin{lemme} Soit ${\mathcal  T}\subset X$ un bout de transversale.
Lorsque $k\geq k_{1}({\mathcal  T}) $ est assez grand, alors toute
famille $j$ continue de $1$-jets de sections holomorphes de
$E^{\otimes k}$ d\'efinie le long de ${\mathcal  T}$ se prolonge en une
section holomorphe continue $\sigma:X\rightarrow  E^{\otimes k}$ telle que
\[  |\sigma|_{\infty} \leq C|j|_{\infty, {\mathcal  T}}, \]
o\`u $C$ est une constante universelle.\label{prolongement holomorphe lamine}\end{lemme}

\textit{D\'emonstration.} Quitte \`a consid\'erer des puissances assez grandes de $E$, nous pouvons supposer que 
les rayons des feuilles sont sup\'erieur \`a $r_0$, que la courbure de Ricci de la m\'etrique 
k\"ahl\'erienne induite par la courbure de la m\'etrique $|.|^{\otimes k}$ est major\'ee par $c_0$,
et que le bout de transversale $\mathcal T$ intersecte les feuilles sur des parties $2\delta_0$-s\'epar\'ees. 
Consid\'erons alors un bout de transversale $\mathcal T'$ qui 
contient $\mathcal T$ dans son int\'erieur, et tel que toute 
feuille l'intersecte sur une partie $\delta _0$-s\'epar\'ee.

Commen\c{c}ons par \'etendre la famille de jets $j$ \`a une 
famille continues $j'$ de $1$-jets d\'efinie sur $\mathcal T'$ telle que~: 
d'une part $j'$ s'annule sur le bord de $\mathcal T'$, 
et d'autre part $ |j'|_{\infty, \mathcal T'} \leq |j|_{\infty,\mathcal T}$. D'apr\`es le Lemme 
\ref{convergence des series fuchsiennes}, nous avons 
\[  |J_1 \tilde{\sigma}(j') - j|_{\infty, \mathcal T} \leq D(\delta_0) |j|_{\infty,\mathcal T},\]
et d'apr\`es \ref{continuite des series fuchsiennes laminees}, $\tilde{\sigma}(j)$ est une section continue 
holomorphe le long des feuilles. Ceci est donc le d\'ebut d'un proc\'ed\'e d'approximations successives 
qui nous donnent le r\'esultat. 

\vspace{0.2cm}

Achevons la d\'emonstration du Th\'eor\`eme~\ref{theoreme de Kodaira lamine}. Soient
$x_{1}$ et $x_{2}$ deux points distincts de $X$. D'apr\`es la
Proposition \ref{prolongement holomorphe lamine}, lorsque $k$ est
assez grand, il existe une section holomorphe $\sigma$ de $E^{\otimes k}$
telle que $\sigma(x_{1})=0$ et $\sigma(x_{2})\neq 0$. Par
ailleurs, lorsque $k$ est assez grand, en tout point $y$ de $X$ il
existe des sections holomorphes $\sigma_{0}^{y},\ldots,
\sigma_{n}^{y}$ telles que $|\sigma_{0}^{y}(y)|=1$ et
$f=(\sigma_{1}^{y}/\sigma_{0}^{y},\ldots,
\sigma_{n}^{y}/\sigma_{0}^{y})=z+O(|z|^{2})$, o\`u $z$ est une
coordonn\'ee centr\'ee en $y$. Il existe donc un voisinage autour de
$y$ o\`u $\sigma_{0}^{y}$ ne s'annule pas et $df$ est injective. Un nombre fini
de ses voisinages recouvrent $X$ par compacit\'e. L'application
\[  [\sigma:\sigma_{0}^{y_{1}}:\ldots:\sigma_{n}^{y_{1}}:\ldots:
\sigma_{0}^{y_{r}}:\ldots:\sigma_{n}^{y_{r}}] \]
est une immersion holomorphe et s\'epare les points $x_1$ et $x_2$.

\begin{remarque} En fait nous d\'emontrons que si $k>>1$, alors l'espace des sections 
holomorphes et continues de $E^{\otimes k} $ est 
de dimension infinie, sauf si 
$\mathcal L$ est une vari\'et\'e compacte. 
La situation diff\`ere donc fortement de celle des vari\'et\'es compactes. \end{remarque}

\subsection{Exemple~: les laminations hyperboliques}
Une lamination par surfaces de Riemann est dite \textit{hyperbolique} si toutes ses feuilles sont rev\^etues par le disque unit\'e. 
\'E.~Ghys a construit des fonctions m\'e\-ro\-mor\-phes sur une lamination par surfaces hyperboliques~\cite{Ghys1}, sous l'hypoth\`ese 
qu'il existe une \textit{transversale totale}.\footnote{Une transversale totale est un bout de transversale dont le bord est vide.}
Ici, nous retrouvons ce r\'esultat sans supposer l'existence d'une transversale totale. 
Nous verrons qu'en fait, toute lamination par surfaces de Riemann hyperboliques d'un espace compact 
admet une \textit{multi-transversale totale} (voir Exemple \ref{intersection avec un hyperplan}). 

\begin{corollaire}\label{immersion hyperbolique} Une lamination par surfaces de Riemann 
hyperboliques d'un espace compact est projective. \end{corollaire}

\textit{D\'emonstration.} Soit $g$ la m\'etrique compl\`ete de courbure $-1$ sur le fibr\'e tangent de chaque feuille. 
D'apr\`es le Th\'eor\`eme d'A.~Candel~\cite{Candel1} et d'A.~Verjovsky~\cite{Ve}, $g$ est lisse sur l'espace total 
de la lamination. Ainsi, le fibr\'e cotangent est muni d'une m\'etrique lisse de courbure strictement positive. 
Le Corollaire r\'esulte donc du Th\'eor\`eme~\ref{theoreme de Kodaira lamine} et du r\'esultat suivant. 

\begin{lemme} Une lamination par surfaces hyperboliques d'un espace compact n'a pas de cycle \'evanouissant. \end{lemme}

\textit{D\'emonstration.} Soit $\mathcal L$ une lamination par surfaces hyperboliques d'un espace compact $X$, 
et supposons que $\mathcal L$
ait un cycle \'evanouissant. Il existe donc un lacet $\gamma_{\infty} : {\bf S}^1 \rightarrow X$ contenu dans une feuille $L$ et non homotope \`a un point
dans $L$, qui est limite uniforme de lacets $\gamma_n :{\bf S}^1 \rightarrow X$ contenus dans des feuilles voisines $L_n$ dans 
lesquelles ils sont homotopes \`a un point. 

Nous pouvons supposer que $\gamma_{\infty}$ est une g\'eod\'esique ferm\'ee. D'autre part, en approximant par des applications lisses 
l'application $\gamma : {\bf S}^1 \times {\bf N}\cup \{\infty\} \rightarrow X$ dans la topologie $C^0$, nous pouvons aussi supposer 
que les courbes $\gamma_n$ sont 
de classe $C^3$ et que la convergence $\gamma_n \rightarrow \gamma$ a lieu dans la topologie 
$C^3$.
Nous en d\'eduisons que la longueur des courbes $\gamma_n$ tend vers la longueur de $\gamma_{\infty}$, et que leur courbure 
tend uniform\'ement vers $0$. 

Le rev\^etement universel des feuilles $L_n$ est le disque hyperbolique ${\bf D}$. 
Comme les $\gamma_n$ sont homotopes \`a un point dans $L_n$, elles 
se rel\`event en des courbes $\widetilde{\gamma_n} :{\bf S}^1\rightarrow {\bf D}$, qui ont m\^eme longueur et m\^eme
courbure que les courbes $\gamma_n$. Comme le disque hyperbolique est homog\`ene, nous pouvons supposer que ces courbes passent 
par un point donn\'e $p$ de $\bf D$. Extrayons une sous-suite de $\widetilde{\gamma_n}$ qui converge dans la topologie $C^2$ 
vers une courbe $\widetilde {\gamma_{\infty}}:{\bf S}^1 \rightarrow {\bf D}$, ce qui est possible car la famille de courbes $\widetilde {\gamma_n}$
est born\'ee dans la topologie $C^3$. Le lacet $\widetilde {\gamma_{\infty}}$ a la m\^eme longueur que $\gamma_{\infty}$, et sa courbure est nulle. 
C'est donc une contradiction, car il n'y a pas de g\'eod\'esique ferm\'ee dans le disque hyperbolique.  

\subsection{Conjecture}\label{conjecture}
\textit{Soit $\mathcal  L$ une lamination
par vari\'et\'es complexes d'un espace compact $X$ dont la dimension topologique~\footnote{Rappelons que la 
dimension topologique d'un espace compact est finie si et seulement si il se plonge 
continument dans un espace euclidien.}
est finie. Supposons qu'il existe un fibr\'e en droites holomorphe
positif $E\rightarrow {\mathcal  L}$, et que $\mathcal L $ n'ait pas de cycle \'evanouissant. 
Alors il existe un \textit{plongement} holomorphe $\pi:X\rightarrow {\bf C}P^N$.} 

\vspace{0.2cm}

T.~Ohsawa et N.~Sibony \cite{Ohsawa-Sibony1} d\'emontrent cette 
conjecture dans le cas des feuilletages lisses 
par vari\'et\'es complexes de codimension $1$ d'une vari\'et\'e compacte. 
Il semble qu'elle soit li\'ee \`a des questions
de r\'egularit\'e transverse des sections holomorphes des puissances de $E$.
Par exemple, en vertu du Lemme suivant, 
nous saurions r\'epondre \`a la conjecture si l'on savait construire des sections holomorphes
\textit{lipshitziennes} passant par une famille de jets donn\'ee
sur un bout de transversale.
\begin{lemme} Soit $T\subset {\bf C}^{r}$ un compact et $B$ la boule unit\'e de
${\bf C}^{n}$. Soit $f:B\times T\rightarrow {\bf C} ^{n+r}$ une application
lipshitzienne holomorphe en $z\in B$, et telle que
\[  f(z,t)=(z,t)+O(|z|^{2}).\]
Alors $f$ plonge un voisinage de $0\times T$.\end{lemme}
\textit{D\'emonstration.} Consid\'erons les coordonn\'ees
\[(Z_{1},\ldots,Z_{n})=(f_{1}(z,t),\ldots,f_{n}(z,t)).\]
L'application $(z,t)\mapsto (Z,t)$ est un biholomorphisme d'un voisinage de
$0\times T$ dans un ouvert de la forme $B(r)\times T$. Dans ces coordonn\'ees
l'application $f$ s'exprime par
\[  f(Z,t)= (Z,t+g(Z,t)),\]
o\`u $g:B(r)\times T\rightarrow {\bf C}^{r}$ est une application holomorphe
et lipshitzienne telle que $g(0,t)=0$ et $dg_{(0,t)}=0$.
Montrons que $f$ est injective au voisinage de $0\times T$.
Le Lemme de Schwarz
appliqu\'e aux fonctions $g(t,.)-g(t',.)$ (dont la norme est major\'ee par $C|t-t'|$)
donne des in\'egalit\'es du type
\[  |g(Z,t)-g(Z,t')|\leq C|t-t'||Z|^{2}.\]
Supposons donc que l'on ait deux points distincts $(Z,t)$ et $(Z',t')$
de $B(r)\times T$ tels que $f(Z,t)=f(Z',t')$. On a bien entendu $Z=Z'$,
d'o\`u l'on d\'eduit l'\'equation $t-t'= g(Z,t')-g(Z,t)$, et les in\'egalit\'es
$|t-t'|\leq  |g(Z,t)-g(Z,t')|\leq C|t-t'||Z|^{2}$. Ceci oblige $|Z|\geq \sqrt{1/C}$
et $f$ plonge $B(\sqrt{1/C})\times T$ dans ${\bf C}^{n+r}$. 

\vspace{0.2cm}

Voici quelques exemples o\`u la m\'ethode des s\'eries fuchsiennes  m\`enent \`a des sections 
lipschitziennes~:

\vspace{0.2cm}

\begin{itemize} 

\item Soit $M$ une vari\'et\'e projective et $\rho : \pi_1(M) \rightarrow T$ une action 
du groupe fondamental de $M$ sur un espace 
m\'etri\-que compact $T$ par transformations bilipschitziennes. 
La \textit{suspension} de $\rho$ est le quotient de la lamination triviale 
$\tilde {M}\times T$ par l'action diagonale du 
groupe fondamental de $M$, o\`u $\tilde {M}$ est le rev\^etement universel de $M$. 
C'est une fibration au dessus de $M$, qui est un rev\^etement holomorphe le long des feuilles.
Soit $E\rightarrow M$ un fibr\'e en droites holomorphe de courbure strictement positive,
et $F\rightarrow M\ltimes _{\rho} T$ le tir\'e en arri\`ere de $E$.  
Il est ais\'e de v\'erifier que 
les s\'eries fuchsiennes construites sur les puissances de $F$ sont lipschitziennes. 
Les suspensions se plongent donc holomorphiquement dans un espace projectif complexe, 
si l'espace $T$ est de dimension topologique finie. 

\item Soit $U$ un domaine sym\'etrique born\'e de ${\bf C}^n$, $G$ son groupe de biholomorphismes, et $K$ 
le stabilisateur d'un point de $U$. Le groupe $G$ est un groupe alg\'ebrique r\'eel qui peut \^etre 
d\'efini sur ${\bf Q}$. Consid\'erons 
un r\'eseau cocompact $\Gamma \subset G({\bf R})\times G({\bf Q_p})$, qui existe d'apr\`es les travaux d'A. Borel 
pour certaines valeurs d'entiers premiers $p$. Observons que $\Gamma$ agit par multiplication \`a gauche sur 
la lamination triviale $U\times G({\bf Q_p}) = G/K \times G({\bf Q_p})$ par biholomorphismes. 
Si l'on suppose que $\Gamma $ agit sans points fixes, alors l'espace quotient est 
compact et muni d'une structure $\mathcal L$ de lamination par vari\'et\'es complexes, dont l'espace transverse est 
totalement discontinu. Le fibr\'e canonique $K$  
de $\mathcal L$ est muni d'une m\'etrique de courbure strictement positive. Les s\'eries fuchsiennes construites 
sur les puissances du fibr\'e canonique m\`enent \`a des sections lipschitziennes. Ces exemples de laminations 
``arithm\'etiques" se plongent donc holomorphiquement dans un espace projectif complexe. 
Nous ne savons pas si ce sont des ensembles limites de feuilletages holomorphes sur des vari\'et\'es projectives. 

\item Sur une lamination par surfaces hyperboliques dont la m\'e\-tri\-que hyperbolique
est \textit{lipschitzienne}, les s\'eries fuchsiennes construites \`a partir des puissances du fibr\'e 
canonique sont lipschitziennes. Ces laminations se plongent donc
holomorphiquement dans un espace projectif complexe. Des exemples de telles laminations apparaissent 
avec les laminations associ\'ees \`a des pavages de l'espace hyperbolique (voir \cite{Ghys1}). 
\`A nouveau, nous ne savons pas si ces laminations sont des minimaux de feuilletages holomorphes 
sur des vari\'et\'es projectives. 

\end{itemize}

\subsection{Plongement symplectique}\label{Plongement symplectique}
L'espace projectif complexe ${\bf C}P^N$ est muni de la $(1,1)$-forme de Fubini-Study
qui est d\'efinie dans les coordonn\'ees homog\`enes $[x_{0}:\ldots:x_{N}]$ par
\[  \omega_{FS} = \frac{\sqrt{-1}}{2\pi}\overline{\partial}\partial
\log (|x_{0}|^{2}+\ldots +|x_{N}|^{2}).\]
C'est une forme \textit{symplectique} compatible avec la structure complexe standard.

Dans ce paragraphe nous d\'emontrons le Th\'eor\`eme \ref{immersion symplectique}. 
Soit $\mathcal  L$ une lamination par vari\'et\'es complexes
d'un espace compact $X$, sans cycle \'evanouissant, de dimension topologique finie et $E\rightarrow {\mathcal  L}$ un
fibr\'e en droites positif. Il s'agit de construire un plongement lisse
$\pi: X\rightarrow {\bf C}P^N$ qui immerge symplectiquement chaque
feuille.

Dans un premier temps, nous
construisons pour tout point $x$ de $X$ une application lisse $\pi
:X\rightarrow {\bf C}P^N$ immergeant symplectiquement chaque
feuille et plongeant un voisinage de $x$. Soit $B\times T'$ une
bo\^{\i}te locale en $x$ pour laquelle $x=(0,t_{0})$ et $T\subset T'$
est un voisinage compact de $t_{0}$ dans $T$. Comme la dimension
topologique de $T$ est finie, il existe un plongement topologique
$p:T\rightarrow {\bf C}^{r}$, pour $r$ assez grand. Si $k$ est
assez grand, il existe des sections holomorphes
$\sigma_{0},\ldots,\sigma_{r+n}$ de $E^{k}$, born\'ees par une
constante $C$ ne d\'ependant que de $T$ telles que
\begin{itemize}
\item $|\sigma_{0}|=1$,
\item $(\sigma_{1}(z,t)/\sigma_{0}(z,t),\ldots,\sigma_{r+n}(0,z)/\sigma_{0}(0,z))
=(p(t),z)+O(|z|^{2})$.
\end{itemize}
On compl\`ete cette famille par des sections holomorphes
$\sigma_{n+r+1},\ldots, \sigma_{N}$ de $E^{k}$ born\'ees par $C$ et telles que
$\pi=[\sigma_{0}:\ldots:\sigma_{N}]$ induise une application continue de $X$ dans
${\bf C}P^N$ immergeant holomorphiquement chaque feuille (voir \ref{demonstration du theoreme de Kodaira lamine}).
Nous allons d\'eformer $\pi$ dans la topologie $C^{1}$ en une application
qui plonge un voisinage de $x$. Comme $\pi$ est symplectique,
nous aurons encore une application symplectique.

Soit $\psi:B\rightarrow {\bf R}$ une fonction lisse, positive,
valant $1$ sur $B(1/3)$ et nulle \`a l'ext\'erieur de $B(2/3)$.
Pour $0<\epsilon<1$ posons $\psi_{\epsilon}(x)=\psi(x/\epsilon)$
pour $x\in B$. Le support de $\psi_{\epsilon}$ est contenu dans
$B(\epsilon/3)$. Pour $\epsilon=0$ on pose $\psi_{0}=0$. Soit
$\epsilon:T\rightarrow {\bf R}$ une fonction positive qui
s'annule sur le bord de $T$ et qui est strictement positive et
$t_{0}$. Ecrivons $p=(p_{1},\ldots,p_{r})$. Soit
$\sigma_{0}'=\sigma_{0}$. Si $1\leq i\leq r$, posons
\[  \sigma_{i}'(z,t)= \psi_{\epsilon}(t)(z)\sigma_{0}(z,t)p_{i}(t)+
(1-\psi_{\epsilon(t)}(z))\sigma_{i}(z,t).\]
Si $1\leq i \leq n$, posons
\[  \sigma_{r+i}(z,t)= \psi_{\epsilon}(t)(z)\sigma_{0}(z,t) z_{i}+
(1-\psi_{\epsilon(t)}(z))\sigma_{r+i}(z,t).\]
Si $r+n\leq i\leq N$, nous posons $\sigma_{i}'=\sigma_{i}$.

Soit $r>0$ un r\'eel tel que $|\sigma_{0}|\geq 1/2$ sur $B(r)$.
Lorsque $|\epsilon|_{\infty}<r$, les sections $\sigma_{i}'$ se
prolongent en des sections lisses de $E^{k}$ par
$\sigma_{i}'=\sigma_{i}$ \`a l'ext\'erieur de $B\times T$.
L'application $[\sigma_{0}':\ldots:\sigma_{N}']$ induit alors une
application lisse de $X$ dans ${\bf C}P^N$, en plongeant un
voisinage de $x$.

Comparons les fonctions
\[f'=(\sigma_{1}'/\sigma_{0}',\ldots, \sigma_{N}'/\sigma_{0}')\ \ \ 
\mathrm{et}\ \ \ f=(\sigma_{1}/\sigma_{0},\ldots, \sigma_{N}/\sigma_{0})\]
dans la topologie $C^{1}$. Elles sont d\'efinies sur $B(r)$
et \`a valeurs dans ${\bf C}^{N}$. Nous avons
\[ f'(z,t)-f(z,t)=\psi_{\epsilon(t)}(z)
(g_{1}(z,t),\ldots,g_{N}(z,t)),\]
avec $g_{i}(z,t)=f_{i}(z,t)-f_{1}(z,t)$ si $1\leq 1\leq r$,
$g_{i}(z,t)=f_{i}(z,t)-z_{i-r}$ si $r+1\leq i\leq r+n$ et $g_{i}(z,t)=0$ sinon.
Les fonctions $g_{i}$ sont des fonctions holomorphes d\'efinies sur $B(r)$,
born\'ees par une constante $C$ et qui s'annule en $0$ ainsi que leur diff\'erentielle.
Nous avons donc des majorations $|g(z,t)|\leq C |z|^{2}$ et
$|dg_{(z,t)}|\leq C|z|$, d'apr\`es le Lemme de Schwarz.
Puisque le support de $\psi_{\epsilon}$ est contenu
dans la boule $B(\epsilon/3)$ nous avons $|f-f'|_{\infty}\leq C|\epsilon|_{\infty}^{2}$.
D'autre part, en un point $t$ tel que $\epsilon (t)>0$ nous avons
\[ df'-df= d\psi_{\epsilon(t)} g+\psi_{\epsilon(t)}dg.\]
Or $|d\psi_{\epsilon}|\leq C/\epsilon$, et donc puisque \`a nouveau le support de
$\psi_{\epsilon}$ est contenu dans $B(\epsilon/3)$
\[ |df'-df|\leq C|\epsilon|_{\infty}.\]
Ceci ach\`eve la d\'emonstration de la premi\`ere \'etape.

Le plongement de Pl\"ucker est un plongement holomorphe et symplectique de 
${\bf C}P^{N_1}\times {\bf C}P^{N_{2}}$ dans ${\bf C}P^{N_3}$, avec 
$N_3=(N_{1}+1)(N_{2}+1)-1$.
Il est d\'efini par la formule
$P:([x_{i}],[y_{j}])\mapsto [x_{i}y_{j}]$. Supposons que l'on ait
deux applications symplectiques $\pi_{i}:X\rightarrow 
{\bf C}P^{N_i}$. Il est imm\'ediat de voir que l'application $P\circ
(\pi_{1},\pi_{2})$ est encore symplectique. En utilisant les
applications de la premi\`ere \'etape, et en faisant des produits
compos\'es par des plongements de Pl\"ucker, on d\'emontre qu'il
existe une immersion lisse $\pi:X\rightarrow
{\bf C}P^N$ symplectique le long des feuilles.

Pour construire un \textit{plongement}, remarquons que l'ensemble 
\[ F=\{(x,y)\in X\times X,\ \pi(x)=\pi (y)\ et \ x\neq y\}\] 
est compact car $\pi$ est une immersion. Le Th\'eor\`eme~\ref{theoreme de Kodaira lamine} nous dit 
que les applications holomorphes \`a valeurs dans ${\bf C}P^N$ immergeant
holomorphiquement chaque feuilles s\'eparent les points de $X$.
Si un couple $(x,y)$ appartient \`a $F$, il existe une application
holomorphe immergeant holomorphiquement chaque feuille
$\pi':X\rightarrow {\bf C}P^{N'}$ telle que
$\pi'(x)\neq \pi'(y)$. On a $\pi'(x')\neq\pi'(y')$ dans un voisinage de $(x,y)$.
En prenant un nombre fini $\pi_{1}',\ldots, \pi_{s}'$, on construit alors
une application $P\circ (\pi_{1},\ldots,\pi_{s})$ plongeant symplectiquement
$X$ dans un projectif complexe.

\vspace{0.2cm}

Soit $(X,{\mathcal  L},\omega)$ une lamination compacte
symplectique, c'est \`a dire que les feuilles de ${\mathcal  L}$ sont
de dimension paire et que $\omega$ est une forme lisse qui est
symplectique le long des feuilles. Supposons qu'il existe un
fibr\'e en cercles au dessus de ${\mathcal  L}$ qui ait une connexion
dont la courbure est $-i\omega$, et que $X$ soit de dimension topologique finie. Alors
nous conjecturons qu'il existe un plongement symplectique $\pi: X\rightarrow {\bf C}P^N$.
Le Th\'eor\`eme \ref{immersion symplectique} d\'emontre
cette conjecture dans le cas d'une lamination par surfaces. 
A. Ibort et D. Mart\^{\i}nez Torres l'ont d\'emontr\'e dans 
le cas d'un feuilletage de codimension $1$~\cite{I-M}. 

\section{Le cas de la dimension $1$}

Nous avons vu qu'une surface riemannienne \`a g\'eom\'etrie born\'ee est projective (Th\'eor\`eme~\ref{surfaces orientees}). 
Cependant, il y a des exemples de laminations par surfaces de Riemann d'un espace compact 
qui ne le sont pas. Dans cette partie, 
nous donnons deux conditions n\'ecessaires et suffisantes pour qu'une lamination par surfaces de Riemann d'un espace compact
(qui n'a pas de cycle \'evanouissant) soit projective. Ces conditions portent sur 
la topologie du feuilletage, et ne tiennent pas compte de
la structure conforme le long des feuilles. 

\subsection{Diviseurs} 
Soit $\mathcal L$ une lamination par surfaces de Riemann d'un espace compact~$X$. 

\begin{definition}
Un \textit{diviseur} est une partie ferm\'ee ${\mathcal D}$ de $X$, avec une fonction 
$m:{\mathcal D}\rightarrow {\bf N} ^{*}$ telle que au voisinage ${\mathcal V}$ de tout point de ${\mathcal D}$ 
il existe une fonction holomorphe $f:{\mathcal V}\rightarrow {\bf C}$ ne s'annulant identiquement sur aucune 
feuille et ayant la propri\'et\'e que ${\mathcal D}\cap {\mathcal V}=f^{-1}(0)$ et que la multiplicit\'e 
d'annulation $m_{t}(f)$ de $f$ en tout point $t$ de ${\mathcal D}$ est \'egale \`a $m(t)$.\end{definition}

\begin{exemple}
\label{intersection avec un hyperplan}
L'intersection d'une lamination par surfaces de Riemann projective avec 
un hyperplan complexe est un diviseur, si la lamination n'a aucune feuille contenue dans l'hyperplan.
Par cons\'equent, sur une lamination par surfaces de Riemann projective, il existe des diviseurs passant par tout point.   

Un exemple de lamination par surfaces qui admet une feuille n'intersectant pas de diviseur 
est la \textit{composante de Reeb}. 
C'est le quotient du feuilletage horizontal de ${\bf C}\times [0,\infty) -\{0\}$ 
par l'homoth\'etie de rapport $2$. Le lecteur pourra s'assurer 
qu'il n'y a pas de diviseur passant par la feuille compacte 
quotient de ${\bf C} \times \{0\} -\{0\}$.\end{exemple} 

Comme dans le cas d'une surface de Riemann compacte, 
\`a un diviseur ${\mathcal D}$ est associ\'e 
un fibr\'e en droites holomorphe $E_{{\mathcal D}}\rightarrow {\mathcal L}$. 
Nous rappelons cette construction ci-apr\`es. 

Soit $\{U_{i}\}$ un 
recouvrement de $X$ par des bo\^{\i}tes sur lesquelles sont donn\'es des \'equations $f_{i}=0$
d\'efinissant ${\mathcal D}$ et $m$, o\`u les $f_{i}:U_{i}\rightarrow {\bf C}$ sont des fonctions holomorphes. 
Nous supposerons dans toute la suite que les fonctions $f_{i}$ ne s'annulent pas sur le bord des plaques 
de $U_{i}$. Les quotients $f_{j}/f_{i}$ d\'efinissent des fonctions holomorphes non nulles $g_{i,j}$ sur chaque
plaque de l'intersection $U_{i}\cap U_{j}$. La formule de Cauchy 
\[  g_{i,j}(z)= \frac{1}{2\pi \sqrt{-1}} \int_{\partial D}\frac{g_{i,j}(\zeta)}{z-\zeta}d\zeta  \]
montre que les $g_{i,j}$ sont en fait continues m\^eme aux points o\`u $f_{i}$ s'annule. 
Bien entendu, nous avons  
les relations de cocycle 
\[  g_{i,j}g_{j,k}=g_{i,k} , \]
qui nous permettent de construire un fibr\'e en droites holomorphe $E\rightarrow {\mathcal L}$ : sur chaque $U_{i}$, 
il y a une section holomorphe $s_{i}:U_{i}\rightarrow E$ ne s'annulant pas et entre deux de ces sections il y a 
la relation 
\[ s_{i}=g_{i,j}s_{j} . \]
Remarquons que l'on a alors $f_{i}s_{i}=f_{j}s_{j}$, si bien qu'il y a une section globale $s$ de $E$ d\'efinie 
sur $U_{i}$ par $f_{i}s_{i}$. Elle s'annule donc exactement sur ${\mathcal D}$ avec la multiplicit\'e $m$ le long des feuilles.

\begin{proposition} Soit ${\mathcal L}$ une lamination par surfaces de Riemann d'un espace compact $X$ et ${\mathcal D}$ 
un diviseur. Alors le fibr\'e en droites $E_{\mathcal D}$ est muni d'une m\'etrique hermitienne 
dont la courbure est positive et non nulle sur ${\mathcal D}$.\end{proposition} 
\textit{D\'emonstration.} Nous devons construire la norme $|s|$ de la section $s$. 
En fait nous construisons plut\^ot la fonction $\varphi =\log |s|$.
Ce doit \^etre une fonction lisse en dehors de ${\mathcal D}$, et elle doit pr\'esenter des singularit\'es 
logarithmiques le long de ${\mathcal D}$ : dans chaque bo\^{\i}te $U_{i}$ les fonctions 
\[ \varphi -\log|f_{i}|\]
sont lisses. Comme on veut que la courbure soit positive, il faut de plus que le laplacien de $\varphi$ soit 
n\'egatif, et strictement dans un voisinage de ${\mathcal D}$.

Soit $\{U_{i}\}$ un recouvrement fini de $X$ par des bo\^{\i}tes 
${\bf D} \times T_{i}$ 
o\`u il y a des fonctions holomorphes $f_{i}:U_{i}\rightarrow {\bf C}$ dont le lieu des
$0$ est ${\mathcal D}$ et dont les multiplicit\'es le long des feuilles est $m$. Nous supposerons que les $f_{i}$
sont born\'ees et que $|f_{i}|\geq \alpha$ \`a l'ext\'erieur de ${\bf D}_{1/2}\times T_{i}$, o\`u $\alpha$ est une 
constante strictement positive.

Sur chaque bo\^{\i}te $U_{i}$, consid\'erons les fonctions 
\[  \chi = \inf(C, -\log |f_{i}|) ,\]
o\`u $C$ est une constante r\'eelle que l'on choisira assez grande. Ce sont des fonctions surharmoniques,
et sont harmoniques sur l'ouvert $V_{C}=\{-\log |f_{i}|<C\} $. Lorsque $C$ tend vers $+\infty$, la famille
$V_{C}$ est une base de voisinage de ${\mathcal D}$. Prenons un noyau r\'egularisant $K$ et d\'efinissons 
\[ \psi _{i}(x) = \int_{{\bf D}} K(x-y)\chi(y)dyd\overline{y}. \]
On demande que $K:{\bf D}\rightarrow {\bf R}$ soit une fonction lisse positive, strictement positive sur 
${\bf D}_{1/3}$ et nulle \`a l'ext\'erieur, et que son int\'egrale soit $1$. Les fonctions $\psi _{i}$ 
sont alors lisses, d\'efinie sur ${\bf D}_{2/3}$ et on a 
\[  \triangle \psi_{i}(x)=\int_{{\bf D}} K(x-y)\triangle\chi dyd\overline{y}(y),\]
o\`u il faut voir $\mu=\triangle\chi dyd\overline{y}$ comme une mesure n\'egative sur ${\bf D}$. Ainsi les 
$\psi_{i}$ sont surharmoniques et on a 
\[  \triangle\psi_{i} <0\]
aux points distants du support de $\mu$ de moins de $1/3$. Pour $C$ assez grand, $\triangle\psi_{i}$ est 
donc strictement n\'egative sur ${\mathcal D}$. Par ailleurs, si $K(x-y)$ ne d\'epend que de la distance  
entre $x$ et $y$, nous avons 
\[ \psi_{i}=\chi=-\log |f_{i}|\]
aux points situ\'es \`a une distance sup\'erieure \`a $1/3$ de $V_{C}$, 
par harmonicit\'e de $\chi$ en dehors de $V_{C}$. Remarquons que pour $C>-\log \alpha$, l'ouvert $V_{C}$ 
est inclu dans ${\bf D}_{1/2}$. Dans ce cas les fonctions $\psi_{i}$ se prolonge \`a tout $U_{i}$ en des 
fonctions lisses valant $-\log |f_{i}|$ \`a l'ext\'erieur de $D_{1/2+1/3}$.

Soit ${\mathcal D}_{i}$ l'intersection de ${\mathcal D}$ et de $U_{i}$. Consid\'erons une partition de l'unit\'e 
$\{\rho_{i}\}$ associ\'ee au recouvrement $\{{\mathcal D}_{i}\}$ de ${\mathcal D}$. Les fonctions 
\[  \varphi_{i}= \rho_{i}(\psi_{i} +\log |f_{i}|)\]
sont des fonctions lisses sur $X-{\mathcal D}$ surharmonique le long des feuilles et pr\'esentent un p\^ole 
\[  \rho_{i} \log |f_{i}|  \]
le long de ${\mathcal D}$. On pose $\varphi =\sum_{i}  \varphi_{i}$ : dans une bo\^{\i}te $U_{j}$ nous avons 
\[ \varphi = \log |f_{j}| + \sum_{i} \rho_{i}(\psi_{i} + \log|g_{i,j}|) . \]
Ceci d\'emontre que $\varphi$ a le p\^ole $\log|f_{i}|$ le long de ${\mathcal D}$. De plus on a 
\[ \triangle \varphi =\sum_{i} \rho_{i} \triangle \psi_{i},  \]
ce qui prouve la surharmonicit\'e de $\varphi$, stricte sur ${\mathcal D}$. La Proposition est d\'emontr\'ee. 

\begin{lemme} Une lamination compacte par surfaces de Riemann sans cycle \'evanouissant 
poss\`ede un fibr\'e en droites strictement positif 
partout si et seulement si il existe un diviseur coupant toutes les 
feuilles. \label{multi-transversales}\end{lemme}

\textit{D\'emonstration.} Nous avons montr\'e que si il y a un fibr\'e 
en droites strictement positif au dessus d'une lamination par surfaces de Riemann compacte $(X,{\mathcal L})$ 
n'ayant pas de cycle \'evanouissant, 
alors il existe une application holomorphe $X\rightarrow {\bf C}P^1$ identiquement constante sur aucune feuille
(voir Th\'eor\`eme~\ref{demonstration du theoreme de Kodaira lamine}). 
Les fibres de cette application sont des diviseurs. La r\'eunion 
d'un nombre fini d'entre eux coupe toute les feuilles. Mais c'est encore un diviseur.

R\'eciproquement, supposons qu'il existe un diviseur ${\mathcal D}$ coupant toutes les feuilles. 
En le d\'epla\c{c}ant le long du feuilletage, nous obtenons 
un diviseur passant par un point donn\'e $x$. En vertu du Lemme 
pr\'ec\'edent, il existe un fibr\'e 
en droites holomorphe hermitien dont la courbure est positive partout et strictement positive en $x$. Elle est donc strictement 
positive sur un voisinage $V_{x}$ de $x$. Un nombre fini de ces voisinage recouvrent $X$~: 
le produit des fibr\'es en droites correspondant 
est un fibr\'e en droites holomorphe hermitien dont la courbure est partout strictement positive. 
Le Lemme est d\'emontr\'e.

\vspace{0.2cm}

Le Th\'eor\`eme \ref{obstruction 1} r\'esulte du 
Th\'eor\`e\-me~\ref{theoreme de Kodaira lamine} et du Lemme~\ref{multi-transversales}.

\subsection{Cycles feuillet\'es}\label{cycle feuillete}

Un \textit{cycle feuillet\'e} est l'analogue feuillet\'e de la classe 
fondamentale d'une vari\'et\'e compacte. C'est une notion d\^ue \`a 
S.~Schwartzman \cite{Schwartzman1} dans le cas des feuilletages de dimension
r\'eelle $1$. Elle a \'et\'e \'etendue par D.~Ruelle et D.~Sullivan \cite{Ruelle-Sullivan1} 
au cas des feuilletages de dimension sup\'erieure, et par D.~Sullivan \cite{Sullivan1} 
\`a d'autres situations dynamiques. 

\begin{definition} Soit ${\mathcal L}$ une lamination compacte lisse orient\'ee de 
dimension $n$ d'un espace compact $X$. Un \textit{cycle feuillet\'e} est un op\'erateur 
$T:\Omega^{n}({\mathcal L})\rightarrow {\bf R}$ strictement positif sur les formes 
strictement positives et nul sur les formes exactes.\end{definition}

Un cycle feuillet\'e poss\`ede une classe d'homologie d\'efinie de la mani\`ere suivante. 
Observons d'abord que si $\pi : X\rightarrow M$ est une application lisse \`a valeurs dans une vari\'et\'e 
$M$ lisse et \`a bord, le courant $\pi_* T$ est un courant ferm\'e de dimension $n$ sur $M$ et sa classe d'homologie  
$T(M,\pi):=[\pi_{*}T] \in H_{n}(M,{\bf R})$ est bien d\'efinie. 

\begin{proposition} Pour tout cycle feuillet\'e $T$, il existe une unique classe d'homologie $[T]$
dans le $n$-i\`eme groupe d'homologie de Cech telle que 
$T(M,\pi)= \pi_* [T]$ pour toute fonction lisse $\pi: X\rightarrow M$ \`a valeurs dans une vari\'et\'e lisse. \end{proposition}

\textit{D\'emonstration.} Bien entendu, si $f:M\rightarrow N$ est une application lisse, nous avons la relation 
\[   f_{*}T(M,\pi)=T(N,f\circ \pi).\]
La classe d'homologie d'un cycle feuillet\'e est donc un \'el\'ement de
la \textit{limite projective} ${\hat H}_{n}({\mathcal L},{\bf R})$ des groupes 
$H_{n}(M,{\bf R})$ lorsque $(M,\pi)$ d\'ecrit l'ensemble des applications lisses $\pi :X\rightarrow M$
\`a valeurs dans une vari\'et\'e ferm\'ee lisse $M$ \`a bord. C'est le sous-ensemble de 
\[  \prod _{(M,\pi)} H_{n}(M,{\bf R}) \]
form\'e des \'el\'ements $(x(M,\pi))_{(M,\pi)}$ tels que $f_{*}x(M,\pi)=x(N,f\circ \pi)$ pour toute 
application lisse $f:M\rightarrow N$. Il est facile de s'apercevoir que le groupe ${\hat H}_{n}({\mathcal L},{\bf R})$ 
est naturellement isomorphe au 
$n$-i\`eme groupe d'homologie de Cech $H_n^{Cech} (X,{\bf R})$ de $X$~; ceci utilise les propri\'et\'es de continuit\'e 
de la cohomologie de Cech (voir \cite{Spanier}), 
et le fait que l'on peut approcher dans la topologie $C^0$ une application continue $\pi: X\rightarrow M$ par 
des applications lisses. La proposition est d\'emontr\'ee. 

\begin{definition}
Une lamination d'un espace compact est dite \textit{tendue} si aucun cycle feuillet\'e n'est homologue 
\`a $0$. \end{definition}

\begin{exemple} [Ghys] 
Une lamination par surfaces de Riemann projective est tendue. 
En effet, soit $\omega $ la forme de Fubini-Study sur ${\bf C}P^N$. C'est une forme ferm\'ee qui est strictement 
positive sur chaque droite complexe du fibr\'e tangent de ${\bf C}P^N$. La forme $\pi ^{*} \omega$ est alors strictement 
positive sur ${\mathcal L}$. Nous avons donc $T(\pi^{*}\omega)>0$ pour tout 
cycle feuillet\'e $T$. Mais $T(\pi^{*}\omega)=\pi_{*}T(\omega)$, ce qui montre que $\pi_{*}T$ est non homologue 
\`a $0$ puisque $\omega$ est ferm\'ee. Ainsi aucun cycle feuillet\'e $T$ de ${\mathcal L}$ n'est homologue \`a $0$.
\end{exemple}

Pour d\'emontrer le Th\'eor\`eme \ref{obstruction 2}, 
nous construisons un fibr\'e en droites holomorphe positif sur toute lamination 
par surfaces de Riemann tendue d'un espace compact. 
Il suffit alors d'appliquer le Th\'eor\`e\-me~\ref{theoreme de Kodaira lamine}.
Nous commen\c{c}ons par construire des fibr\'es 
en droites complexes hermitiens au dessus de vari\'et\'es lisses, ayant une connexion compatible avec la m\'etrique 
et dont la courbure est une $2$-forme ferm\'ee enti\`ere donn\'ee. 
Ceci est bien connu dans le cas des vari\'et\'es, mais nous le rappelons pour pouvoir l'adapter au cas d'une lamination 
par surfaces. 

\subsubsection{Fibr\'es en cercles}\label{Fibre en cercles} Soit $M$ une vari\'et\'e lisse, compacte et \`a bord. 
Un fibr\'e en cercles lisse au dessus de $M$ est 
une fibration lisse ${\bf S} ^{1}\rightarrow F\rightarrow M$ de groupe structural ${\bf S}^{1}$. 
Soit ${\mathcal U}=\{U_{i}\}$ un recouvrement de $M$ par des ouverts sur lesquels sont d\'efinies des sections locales 
\[ s_{i}:U_{i}\rightarrow F . \]
Nous pouvons red\'ecouvrir la fibration $F\rightarrow M$ par le cocycle de fonctions 
$f_{i,j}:U_{i}\cap U_{j}\rightarrow {\bf S}^{1}$ d\'efini par les relations
\[ s_{i}=f_{i,j}s_{j} \]
sur l'intersection $U_{i}\cap U_{j}$. 


Une connexion lisse sur $F$ est une mani\`ere d'identifier des fibres infinit\'esimalement proches de $F$. C'est la 
donn\'ee d'un op\'erateur $\nabla$ qui \`a une section locale lisse $s$ associe un objet de la forme 
\[  A.s ,\]
o\`u $A$ est une $1$-forme \`a valeurs imaginaires pures, et v\'erifiant des relations du type 
\[ \nabla (fs)= df.s +f \nabla s,\] 
pour toute fonction lisse $f$ et toute section lisse $s$. Posons 
\[  \nabla s_{i}= A_{i} s_{i}.  \]
Nous avons 
\[ \nabla s_{i}= df_{i,j} s_{j} +f_{i,j} A_{j} s_{j}= (d\log f_{i,j} +A_{j})s_{i}, \]
d'o\`u la relation 
\[  A_{i}= d\log f_{i,j}+A_{j}.  \]
La courbure de la connexion $\nabla $ est une $2$-forme ferm\'ee \`a valeurs imaginaires pures d\'efinie par la formule 
\[ \Omega = d A_{i} . \]
Lorsque $\Omega$ est la forme de courbure d'une connexion d'un fibr\'e en cercles, 
on d\'emontre facilement que la classe de cohomologie $\frac{\sqrt{-1}}{2\pi} [\Omega]$ est enti\`ere, 
c'est \`a dire qu'on obtient un entier en int\'egrant $\frac{\sqrt{-1}}{2\pi} \Omega $ sur une 
surface compacte immerg\'ee dans $M$. Voici la r\'eciproque de ce r\'esultat.
\begin{lemme}\label{construction d'un fibre dans une classe de cohomologie} 
Soit $\Omega $ une $2$-forme lisse \`a valeurs imaginaires pures sur $M$ telle que 
$\frac{\sqrt{-1}}{2\pi} [\Omega]$ soient enti\`ere. Alors il existe un fibr\'e en cercles lisse 
${\bf S} ^{1}\rightarrow F\rightarrow M$ avec une connexion $\nabla $ dont la courbure est la forme $\Omega$.
\end{lemme}
\textit{D\'emonstration.} Il existe un recouvrement ${\mathcal U}=\{U_{i}\}$ de $M$ par des ouverts $U_{i}$ o\`u 
la forme $\Omega $ est exacte, c'est \`a dire que 
\[ \Omega = dA_{i}  \]
pour une certaine $1$-forme $A_{i}$ \`a valeurs imaginaires pures d\'efinie sur $U_{i}$. Dans ce cas les $1$-formes 
$A_{i}-A_{j}$ sont ferm\'ees sur l'intersection $U_{i}\cap U_{j}$, et quitte \`a r\'eduire les $U_{i}$, on peut les 
supposer exactes. C'est donc qu'il existe des fonctions \`a valeurs imaginaires pures qui sur l'intersection 
$U_{i}\cap U_{j}$ v\'erifient 
\[ d\varphi _{i,j}= A_{i}-A_{j} .\]
Sur l'intersection $U_{i}\cap U_{j}\cap U_{k}$ de trois des ouverts du recouvrement, que nous pouvons supposer connexe,
les fonctions  
\begin{equation}\label{cocycle} \varphi _{i,j}+\varphi _{j,k} +\varphi _{k,i}  \end{equation}
sont constantes. On d\'efinit ainsi un cocyle de $H^{2}(M,{\bf R})$ qui repr\'esente la classe de cohomologie de 
$\Omega$. Comme $\frac{\sqrt{-1}}{2\pi} [\Omega]$ est enti\`ere, c'est que, quitte \`a ajouter 
aux $\varphi_{i,j}$ des constantes, les valeurs des fonctions $(\ref{cocycle})$ sont des multiples entiers de $2\pi \sqrt{-1}$. 
Posons alors 
\[ f_{i,j}= e^{\varphi_{i,j}}.  \]
Ce sont des fonctions \`a valeurs dans le cercle et nous avons la relation de cocycle 
\[  f_{i,j}f_{j,k}f_{k,i}=1 .\]
Ceci nous permet de construire un fibr\'e en cercle $F\rightarrow M$ localement trivial sur les $U_{i}$ avec des sections
trivialisantes $s_{i}$ v\'erifiant 
\[ s_{i}=f_{i,j}s_{j} .\]
D\'efinissons alors une connexion $\nabla$ par les formules 
\[ ds_{i}= A_{i}s_{i}.\]
Sa courbure est la forme $\Omega $. Le Lemme est d\'emontr\'e. 
\subsubsection{Fibr\'e en droites holomorphes} A partir d'un fibr\'e en cercles lisse ${\bf S}^{1}\rightarrow F\rightarrow M$, 
nous construisons naturellement 
un fibr\'e en \textit{droites complexes} lisse ${\bf C}\rightarrow E\rightarrow M$, avec une m\'etrique hermitienne 
$|.|$ lisse dont les sph\`eres unit\'es d\'ecrivent les fibres du fibr\'e en cercles de d\'epart. Les sections de $E$ 
s'\'ecrivent localement 
\[  f_{i} s_{i}, \]
o\`u $f_{i}:U_{i}\rightarrow {\bf C}$ est une fonction \`a valeurs complexes. Si le fibr\'e en cercles $F$ a une 
connexion $\nabla$, elle s'\'etend naturellement \`a une connexion sur $E$ par la formule 
\[ \nabla (fs) = df s+ f \nabla s ,\]   
pour toute fonction lisse $f$ et toute section lisse $s$ de $F$. La m\'etrique hermitienne $|.|$ est alors invariante 
par la connexion. 

Supposons que ${\mathcal L}$ soit une lamination par surfaces de Riemann d'un espace compact $X$, et que l'on ait une application 
\textit{lisse} 
\[ \pi: X\rightarrow M \]
immergeant chaque feuille de ${\mathcal L}$ dans $M$. Le fibr\'e $E\rightarrow M$ se rappatrie en un fibr\'e par droites 
complexes. Ce fibr\'e est lisse et on le note encore $E$. Nous montrons au Lemme~\ref{lemme de Poincare} 
qu'il y a une unique fa\c{c}on de 
le munir d'une structure de fibr\'e en droites \textit{holomorphe} compatible avec la connexion $\nabla$.

Une section locale $s:U\subset X\rightarrow E$ est holomorphe pour une structure compatible avec la connexion 
$\nabla$ si $\nabla s$ est une $(1,0)$-forme \`a valeurs dans $E$ le long de ${\mathcal L}$. 
Supposons que $s_{1}$ et $s_{2}$ soient deux sections holomorphes, avec $s_{2}$ non nulle, et \'ecrivons 
\[ s_{1}= fs_{2},\]
pour une fonction $f$ \`a valeurs dans ${\bf C}$. Nous avons alors 
\[ \nabla s_{1} = df s_{2}+f\nabla s_{2}. \]
Comme les $\nabla s_{i}$, pour $i=1,2$ sont des $(1,0)$-formes, $df$ l'est aussi. C'est donc que $f$ est holomorphe le long des 
feuilles. Ainsi, modulo l'existence locale de sections holomorphes, on d\'efinit sur 
$E\rightarrow {\mathcal L}$ une structure de fibr\'e en droites holomorphe, et c'est de surcroit la seule qui soit compatible 
avec $\nabla$.
\begin{lemme} En tout point de $X$ il existe une section lisse holomorphe ne s'annulant pas.\label{lemme de Poincare}\end{lemme}

\textit{D\'emonstration.} Prenons une section locale $s:U\subset X\rightarrow F$ de $F$, c'est \`a dire une 
section de $E$ de norme $1$. Nous cherchons une fonction $f:U\rightarrow {\bf C}^*$ telle que $\nabla(fs)$ soit une 
$(0,1)$-forme. \'Ecrivons 
\[  \nabla s= A s,\]
pour une $1$-forme lisse $A$ d\'efinie sur $U$. On demande donc que la forme 
\[  d\log f +A \]  
soit de type $(1,0)$. Consid\'erons la partie $(0,1)$ de la forme $A$, que l'on note $A_{a}$. Nous cherchons $f$ en sorte 
que
\begin{equation} \label{inversion du d-barre} \overline{\partial}\log f+ A_{a}=0 .\end{equation}
Nous avons donc \`a inverser le $\overline{\partial}$ avec un param\`etre. Ceci est bien connu : c'est le Lemme 
de $\overline{\partial}$ d\^u \`a Poincar\'e. Si l'on \'ecrit dans une coordonn\'ee holomorphe $z$ la forme $A_{a}$ :
\[ A_{a} = gd\overline{z}, \]
alors une solution au probl\`eme \ref{inversion du d-barre} est donn\'e sous forme d'une int\'egrale par la formule de Poincar\'e :
\[  \log f=\frac{1}{2\pi\sqrt{-1}}\int_{{\bf D}}\frac{g(w)}{w-z}dw\wedge d\overline{w}  .\]   
Le lecteur pourra consulter par exemple : \cite{G-H} page 5. Ceci ach\`eve la d\'emonstration du Lemme.

\subsubsection{Th\'eorie de Sullivan}
Nous adaptons la th\'eorie de Sullivan~\cite{Sullivan1} \`a une lamination compacte orient\'ee plong\'ee 
dans une vari\'et\'e lisse. Le r\'esultat qui nous int\'eresse est le suivant. 

\begin{proposition} Soit $(X,{\mathcal L})$ une lamination compacte, lisse, orien\-t\'ee, de dimension 
topologique finie. Si ${\mathcal L}$ est tendue, il existe un plongement lisse $\pi :X\rightarrow M$ 
\`a valeurs dans une vari\'et\'e lisse, et une $n$-forme lisse sur $M$, ferm\'ee telle que $\pi^*\omega$ est 
strictement positive.\label{approximation} \end{proposition}

Consid\'erons une lamination lisse compacte orient\'ee de dimension 
topologique finie $(X,{\mathcal L})$ et un plongement lisse $\pi: X\rightarrow M$ \`a 
valeurs dans une vari\'et\'e $M$. 
La th\'eorie de Sullivan consiste \`a \'etudier la configuration que forme 
dans l'espace des courants sur $M$:

\vspace{0.2cm}

- L'espace $E$ des courants exacts. 

\vspace{0.2cm}

- L'espace $F$ des courants ferm\'es ($E\subset F$).

\vspace{0.2cm}

- Le c\^one $\mathcal C$ des courants de dimension $n$ sur $M$ 
strictement positifs sur ${\mathcal L}$ (c'est \`a dire que ce sont des courants $T$ sur $M$ 
v\'erifiant la propri\'et\'e suivante: si $\omega$ est une $n$-forme lisse sur $M$ telle 
que $\pi^*\omega$ est strictement positive,
alors $T(\omega)$ est stritement positif). 

\vspace{0.2cm}

Nous adaptons cette th\'eorie pour certains plongements $\pi: X\rightarrow M$~:

\begin{lemme} Il existe un plongement lisse
$\pi: X\rightarrow {\bf R}^N$ v\'erifiant la propri\'et\'e suivante:

\vspace{0.2cm}

(${\mathcal P}$) Pour tout point $x$ de $X$, il existe une coordonn\'ee lisse $(z,t)$ centr\'ee en $x$ 
et une projection lisse $p:{\bf R}^N \rightarrow {\bf R}^n\times {\bf R}^q$
($q=2\times \mathrm{dimension\ topologique\ transverse}+1$) telle que 
\[ p\circ \pi = (z,\tau (t)),\]
o\`u $\tau$ est un plongement d'une transversale locale en $x$ dans ${\bf R}^q$.\end{lemme}

\textit{D\'emonstration.} Dire que $X$ est de dimension topologique finie, c'est exactement dire que
$X$ se plonge topologiquement dans un espace euclidien. Nous voulons rendre ce plongement lisse le long des feuilles.
Donnons nous un point $p$ et une carte locale $U=B \times T$ sur laquelle est d\'efinie une coordonn\'ee
lisse $z$ centr\'ee en $p=(0,t_{0})$.
Choisissons
\begin{itemize}
\item une fonction plateau $\rho: {\bf D}\rightarrow {\bf R} $
lisse, prenant des valeurs entre $0$ et $1$ et telle que $\rho
^{-1}(1) = \{|z|\leq 1/2\}$ et $\rho ^{-1}(1) = \{3/4\leq |z|\}$.
\item une fonction continue $\psi:T\rightarrow {\bf R}$ \`a
valeurs comprises entre $0$ et $1$, telle que $\psi^{-1}(1)$ est un
voisinage de $t_{0}$ et dont le support est un compact inclu dans
$T$. 
\item un plongement $\tau :T\rightarrow {\bf R}^{q}$, o\`u $q=2\times 
\mathrm{dimension\ topologique}(T)+1$.
\end{itemize}
Posons
\[ \pi _{p} (z,t)= \rho (z) \psi (t) (z,\tau (t), 1) \in {\bf R}^{n}\times {\bf R}^{q} \times {\bf R} .\]
C'est une fonction lisse bien d\'efinie sur $X$, et sa restriction \`a $V= \{|z|\leq 1/2\}\times  \psi^{-1}(1)$ est un
plongement. Par compacit\'e, on peut trouver un nombre fini de ces voisinages recouvrant $X$. Le produit des applications
correspondantes $\pi_p$ est alors un plongement lisse de $X$ dans un espace euclidien, qui v\'erife la propri\'et\'e $(\mathcal P)$. 

\begin{lemme} Soit $\mathcal L$ une lamination lisse orient\'ee d'un espace compact $X$ et $\pi:X\rightarrow {\bf R}^N$ 
un plongement lisse v\'erifiant la propri\'et\'e $(\mathcal P)$. 
Un courant ferm\'e agissant sur les $2$-formes lisses \`a support compact de ${\bf R}^N$ et qui est positif sur 
l'image de $X$ est l'image par $\pi$ d'un cycle feuillet\'e 
de $\mathcal L$.\label{courant lamine}\end{lemme}

\textit{D\'emonstration.} Un \'el\'ement de $\mathcal C$ est strictement positif le long de $\mathcal L$ 
et v\'erifie une in\'egalit\'e du type 
\begin{equation}\label{positivite de C} |T(\omega)|\leq D|\omega|_{{\mathcal L},\infty}  ,\end{equation}
pour toute $n$-forme lisse $\omega$ de ${\bf R}^N$, et pour une constante $D$ ne d\'ependant pas de $\omega$. 
Nous voulons d\'efinir $T$ pour \textit{toute forme lisse de $\mathcal L$}:
nous d\'emontrons pour cela le r\'esultat d'approximation suivant.

\vspace{0.2cm}

\noindent \textit{Soit $\pi : X\rightarrow M$ un plongement lisse de $X$ dans une vari\'et\'e lisse $M$, v\'erifiant 
la propri\'et\'e $({\mathcal P})$. Soit $\omega$ une forme lisse sur $\mathcal L$. Alors il existe une 
suite de formes lisses $\omega_n$ sur $M$ telles que $\pi^*\omega_n$ tend vers $\omega$ dans la topologiqe $C^1$.}

\vspace{0.2cm}

\textit{D\'emonstration.} Prenons un recouvrement fini de $X$ par des bo\^{\i}tes $B_i\times T_i$ pour lesquelles il existe 
des coordonn\'ees $(z_i,t_i)$ et des projections lisses $p_i:{\bf R}^N\rightarrow {\bf R}^n\times {\bf R}^q$ 
telles que 
\[ p_i\circ \pi =(z_i,t_i).\] 
En utilisant une partition de l'unit\'e (provenant de fonctions lisses sur $M$), nous pouvons supposer que le support 
de $\omega$ est dans l'une des bo\^{\i}tes $B_i\times T_i$. Pour simplifier les notations nous la noterons $B\times T$ et $(z,t)$ 
les coordonn\'es lisses. Ecrivons alors 
\[\omega =\sum \omega_I dz_I,\]
o\`u les $\omega_I$ sont des fonctions lisses sur $\mathcal L$ (lisses en $z$ dont les d\'eriv\'ees partielles 
par rapport aux variables $z$ sont continues en $z$ et $t$). 


Soit $K$ l'image de $T$ par $\tau$. Nous consid\'erons une famille de mesures de probabilit\'e d\'efinies sur $K$ 
telles que, pour toute fonction continue $f:K\rightarrow {\bf R}$, la fonction 
${\tilde f}:{\bf R}^q\rightarrow {\bf R}$ d\'efinie par 
\[ {\tilde f}(y)=\int_K f(t) d\nu_y(t)\]
soit un prolongement continu de $f$ \`a ${\bf R}^q$. Une telle famille de mesures existe, car il suffit de 
prolonger continument \`a ${\bf R}^q$ l'application $y\in K\mapsto \delta_y \in Prob(K)$, en utilisant des 
approximations simpliciales (pour une d\'emonstration plus constructive, on pourra consulter \cite{St}). 
Ici $Prob(K)$ d\'esigne l'espace des mesures de probabilit\'es sur $K$ et $\delta_y$ est 
la mesure de Dirac en $y$. Nous posons alors 
\[ {\tilde \omega}_I(z,y)=\int_K \omega_I (z,t)d\nu_y(t),\]
pour tout $(z,y)\in B\times {\bf R}^q$. Ce sont des prolongements des $\omega_I$ \`a $B\times {\bf R}^q$, lisses en 
$z$ et dont les d\'eriv\'ees partielles par rapport \`a $z$ sont continues. Il ne reste plus qu'\`a les lisser 
\`a l'aide d'un noyau r\'egularisant $K_{\epsilon}(y,y')$ d\'efini pour $y,y'$ dans ${\bf R}^q$ et $\epsilon$ positif. 
Nous posons 
\[{\tilde \omega}^{\epsilon}_I(z,y)= \int_{{\bf R}^q}K_{\epsilon}(y,y'){\tilde \omega}_I(z,y')d_eucl(y').\]
Il est imm\'ediat de voir que les ${\tilde \omega}^{\epsilon}_I$ sont des fonctions lisses, dont le support est d\'efini 
dans un voisinage aussi petit que l'on veut de $B\times K$. Les formes 
\[ {\tilde \omega}^{\epsilon}=\sum {\tilde \omega}_I ^{\epsilon}dz_I\]
v\'erifient alors 
\[ |\omega -(p\circ \pi)^*{\tilde \omega}^{\epsilon}|_{C^1}\rightarrow 0.\]

\vspace{0.2cm}

Nous sommes en mesure d'achever la d\'emonstration du Lemme~\ref{courant lamine}. D\'efinissons un courant 
$\tilde T$ sur $\mathcal L$ en posant 
\[ {\tilde T}(\omega)= \lim_n T(\omega_n)\]
pour n'importe quelle suite de formes lisses $\omega_n$ sur $M$ telles que $\pi^* \omega_n$ tende vers $\omega$ 
dans la topologie $C^0$. Ceci est bien d\'efini \`a cause de l'in\'egalit\'e \ref{positivite de C}
et du Lemme \ref{approximation}. 
D\'emontrons alors que $\tilde T$ est ferm\'e. Pour cela, prenons une forme exacte $\omega =d\eta$ sur $\mathcal L$ 
et une suite $\eta_n$ de formes lisses sur $M$ telles que $\pi^*\eta_n$ tende vers $\eta$ dans la topologie $C^{1}$
donn\'ee par le Lemme \ref{approximation}. La suite $\omega$ est alors la limite uniforme des formes $\pi^*d\eta_n$. 
Nous avons donc 
\[ {\tilde T}(d\eta) =\lim _n T(d\eta_n)=0.\]
Ceci ach\`eve la d\'emonstration du Lemme. 
 
\vspace{0.2cm}

\textit{D\'emonstration de la Proposition~\ref{approximation}.} 
Prenons un plongement lisse de $X$ dans ${\bf R}^N$ v\'erifiant la propri\'et\'e 
$\mathcal P$. Un cycle feuillet\'e de $\mathcal L$ \'etant non homologue \`a $0$ dans $\mathcal L$, les propri\'et\'es
de continuit\'e de la cohomologie de Cech montrent qu'il existe un voisinage de l'image de $X$ dans ${\bf R}^N$ dans 
lequel il est toujours non homologue \`a $0$. Ceci reste vrai pour des cycles feuillet\'es proche dans la topologie 
faible. Or l'ensemble des cycles feuillet\'es normalis\'es pour la topologie faible est compact. Il existe donc un 
voisinage de l'image de $X$ dans lequel aucun cycle feuillet\'e de $\mathcal L$ n'est homologue \`a $0$. Quitte \`a 
restreindre convenablement ce voisinage on peut supposer que c'est une vari\'et\'e compacte \`a bord lisse. 
Notons la $M$, et $\pi:X\rightarrow M$ le plongement de $X$ dans $M$.

La d\'emonstration est alors identique \`a celle de Sullivan (\cite{Sullivan1}, p. 231). Elle repose sur la dualit\'e 
entre les formes lisses de $M$ et les courants~\cite{Schwartz1}. Remarquons que l'espace $F$ des courants ferm\'es 
est ferm\'e dans l'espace des courants, puisque la diff\'erentiation est continue. De plus, 
l'espace $E$ des courants exacts est ferm\'e dans $F$, puisqu'il est donn\'e par l'annulation des p\'eriodes 
$\omega_i$, les $\omega_i$ formant une base (finie) de $H^i(M,{\bf R})$. D'apr\`es la Proposition 
\ref{courant lamine}, le c\^one $\mathcal C$ des courants ferm\'es et strictement positifs sur l'image de $\mathcal L$
n'intersecte pas $E$. De plus il est \`a base compacte, c'est \`a dire qu'il existe un hyperplan affine qui 
l'intersecte suivant un ensemble compact. Le Th\'eor\`eme de Hahn-Banach montre donc qu'il existe 
une forme lin\'eaire qui est strictement positive sur $\mathcal C$ et nulle sur $E$. Cette forme lin\'eaire 
correspond naturellement \`a une forme lisse, qui est strictement positive sur l'image de $\mathcal L$ et ferm\'ee. 

\vspace{0.2cm}

\textit{D\'emonstration du Th\'eor\`eme \ref{obstruction 2}.} 
Soit ${\mathcal L}$ une lamination tendue par surfaces de Riemann 
d'un espace $X$ de dimension topologique finie. Nous savons qu'il existe 
un plongement lisse 
\[ \pi : X\rightarrow M \]
\`a valeurs dans une vari\'et\'e lisse compacte \`a bord, 
et une $2$-forme $\omega$ lisse sur $M$ strictement positive sur l'image de ${\mathcal L}$ (voir Proposition 3.12). 
La classe de cohomologie dans $H^{2}(M,{\bf R})$ est approximable par des classes de 
cohomologie rationnelles $[\omega_{k}]\in H^{2}(M,{\bf Q})$. On peut supposer que lorsque $k$ tend vers l'infini, 
\[  |\omega -\omega_{k}|_{\infty} \rightarrow 0, \]
quitte \`a choisir $\omega_{k}$ dans sa classe en sorte que $\omega -\omega_{k}$ soit harmonique par rapport \`a une m\'etrique 
riemannienne fix\'ee sur $M$. Pour des entiers $k$ suffisament grand les formes $\omega_{k}$ sont strictement positives sur 
${\mathcal L}$. Un multiple entier d'une d'entre elle est alors une $2$-forme ferm\'ee $\beta$, enti\`ere et strictement positive 
sur ${\mathcal L}$. La Proposition d\'ecoule alors du Lemme~\ref{construction d'un fibre dans une classe de cohomologie}.

\subsection{Crit\`ere pour l'existence d'une multi-transversale totale}

Soit $\mathcal L$ une lamination lisse d'un espace compact $X$. Une \textit{multi-transver\-sa\-le totale} est une partie ferm\'ee 
$\mathcal M \subset X$ et une 
fonction $m : \mathcal M \rightarrow {\bf N}$ telle que au voisinage de tout point de $\mathcal M$ il existe une bo\^{\i}te
$B\times T$, pour laquelle la fonction \[t\mapsto \sum _{(x,t) \in \mathcal M\cap B\times \{t\}} n(x,t)\] est constante.  
Si $\mathcal L$ est une lamination par surfaces de Riemann, le support d'un diviseur est une multi-transversale totale.

Puisque le support d'un diviseur est une multi-transversale totale, le Th\'eor\`eme \ref{theoreme topologique} 
d\'ecoule des Th\'eor\`emes \ref{obstruction 1} et \ref{obstruction 2}. Nous ne savons pas si il est vrai 
en dimension sup\'erieure \`a $3$.

\end{document}